\documentclass[11pt]{amsart}
\usepackage{graphicx}

\usepackage{color}
\usepackage{subfigure}

 \usepackage[
      pdftex=true, 
      colorlinks,
      linkcolor=blue,
      anchorcolor=blue,
      citecolor=blue,
      urlcolor=blue
      ]{hyperref}%

\graphicspath{{./figures/}}

\begin{document}
\newcommand{\bq}{\begin{equation}}
\newcommand{\eq}{\end{equation}}
\newcommand{\e}{\epsilon}
\newcommand{\grad}{\nabla}
\newcommand{\R}{\mathbb{R}}
\newtheorem{theorem}{Theorem}
\newtheorem{prop}{Proposition}
\newtheorem{lemma}[theorem]{Lemma}
\theoremstyle{definition}
\newtheorem{definition}{Definition}
\newtheorem{example}[theorem]{Example}
\newtheorem{xca}[theorem]{Exercise}
\theoremstyle{remark}
\newtheorem{remark}[theorem]{Remark}

%
\newcommand{\x}{\mathbf x}

\newcommand{\C}{\mathcal{C}}
\newcommand{\K}{\mathcal{K}}
\newcommand{\vt}{v_\theta}
\newcommand{\lp}{\lambda_n}
\newcommand{\lm}{\lambda_1}
\newcommand{\dtt}{\frac{d^2 u}{d\vt^2}}

\newcommand{\abs}[1]{\left\vert#1\right\vert}

\newcommand{\dx}{h}
\newcommand{\dth}{d\theta}
\newcommand{\blue}[1]{{\color{blue}{#1}}}


\title[Variational problems with convexity constraints]
{
A numerical method for variational problems with convexity constraints
}
\author{Adam~M. Oberman}
\address{Department of Mathematics, Simon Fraser University}
\email{aoberman@sfu.ca}
\date{\today} 
\begin{abstract}
We consider the problem of approximating the solution
of variational problems subject to the constraint that the admissible
functions must be convex. This problem is at
the interface between convex analysis, convex optimization,
variational problems, and partial differential equation techniques.

The approach is to approximate the (non-polyhedral) cone of convex
functions by a polyhedral cone which can be represented by linear
inequalities.  This approach leads to an optimization problem with linear constraints which can be computed efficiently, hundreds of times faster than existing methods.  
\end{abstract}
   



\maketitle
\tableofcontents

\section{Introduction}

In this article we consider the problem of approximating the solution
of variational problems subject to the constraint that the admissible
functions must be convex. This is a numerical approximation problem at
the interface between convex analysis, convex optimization,
variational problems and Partial Differential Equation (PDE) techniques.

In the theoretical setting, including a convexity constraint
poses no additional difficulties, since the cone of convex functions is itself a convex set.
However, in a computational setting, this problem has proven to be
surprisingly challenging.  First and foremost is the lack of a
computationally tractable characterization of the cone of convex
functions.  Second, there are various mathematical difficulties which
arise when working with approximations of convex functions.  Convexity is not stable under perturbations: while strictly convex functions
are still convex under small perturbations, nonstrictly convex
functions are not.

In this work we present a polyhedral
approximation of the cone of convex functions.  This approximation is
computationally efficient in the sense that it can be represented by a
number of linear inequalities which is small compared to the size of
the problem.  It can be used to build inner (strictly convex) or outer
(slightly nonconvex) approximations to the cone of convex functions.

Our methods are computationally efficient: the computational time to solve the problem with convexity constraints is comparabe to a problem with simpler linear constraints.  The increased efficiency is due to the reduction in the number of constraints to enforce (approximate) convexity.  The results are significantly (hundreds of times) faster, and generally more accurate compared to the method of~\cite{Morin1}, which is currently the most efficient method.

\subsection{Applications and related work}
 The earliest application of variational problems with convexity constraints is Newton's problem of finding a body of minimal resistance~\cite{NewtonsProblem}.  There the convexity constraint
arises as a natural assumption on the shape of the body; see~\cite{NewtonsProblem} for a discussion and also see~\cite{LachPelNewton}.  
A modern application is to mathematical economics~\cite{RochetChone,MMdesign}.  These problems can
often be recast as the projection of a function (in the $L^2(\Omega)$
or $H^1(\Omega)$ norm) onto the set on convex functions defined on the
bounded domain $\Omega$.  Geometric applications include Alexandroff's
problem and Cheeger's problem; see~\cite{LACHAND-ROBERT-05} for
references.  For a discussion of applications to economics and history of this problem, we refer to~\cite{EkelandMB}.

There have been a few different numerical approaches to this problem,
which rely on adapting PDE techniques to the problem at hand.  Early
work~\cite{BrighiChipot} using PDE-type methods, did not make
assertions about the convexity (or approximate convexity) of the
resulting solutions. Later work by~\cite{CLRM}
and~\cite{LACHAND-ROBERT-05} identified some of the difficulties in
working with convex functions.  These difficulties suggest that a
straightforward adaptation of standard numerical methods is not
possible.  The introduction of a large number (superlinear in the
number of variables) of global constraints was required in order to
ensure discrete convexity.

A recent work on the problem is~\cite{Morin1} (see
also~\cite{Morin2}).  In these works, approximate convexity is sought by
enforcing positive definiteness of the discrete Hessian.  However, as the authors of this work explain (see also \S\ref{sec:counter}), the fact that a discrete Hessian is non-negative definite does not ensure that the corresponding points can be interpolated to a convex function.  
The resulting optimization problem
is a conic problem, which is generically more difficult to solve than
one with linear constraints.  However the number of constraints is
less, on the order of the number of variables, in contrast
to~\cite{CLRM}, in which the number of (linear) constraints grew
superlinearly in the number of variables.

In~\cite{HINTERBERGER-06} approximations are performed in the space of bounded Hessians, and as a result, convexity may fail pointwise.

A natural characterization of convex functions on scattered data can
be found in Boyd and Vandenberghe~\cite[\S6.5.5]{BoydBook}: it uses the supporting hyperplane condition.
However it requires the introduction of new variables, and results in
a very large number of constraint equations, one for each pair of data
points.  In~\cite{CLRM}, the variational problem from~\cite{RochetChone} was solved numerically, and convex envelopes were also computed.  
The supporting hyperplane condition for convexity is used to obtain discrete convexity constraints.  This  initially also leads to a quadratic number of constraint equations.   By taking advantage of the fact that points lie on a regular grid, the number of equations is reduced, but even after the reduction, the number of constraint equations is still superlinear in the number of grid points.  

More recently, \cite{EkelandMB} used the supporting hyperplane condition to derive convexity constraints.   In this work the number of constraints is quadratic in the number of points used.  In addition, the function is defined globally, essentially by the Legendre transform, so evaluation of function values and gradients is expensive.   This idea of lifting hyperplanes to satisfy convexity conditions was also used in an early paper~\cite{olikerprussner88} on numerical solutions of the Monge-Amp\`ere partial differential equation.

\section{Counterexamples in approximating convex functions}\label{sec:counter}

The difficulty in working with discrete approximations of convex
functions is that the set, ${\C}$, of convex functions, is a
non-polyhedral cone.  This leads to inequality
constraints, which are more challenging to work with than equality
constraints. In particular, equality constraints have been treated
more frequently in the numerical PDE literature, e.g., the
incompressibility constraint in the Stokes equation.  

The polyhedral convex functions are on the boundary of $\C$:
arbitrarily small perturbations of the these functions can be
nonconvex.  These functions, along with smooth but non-strictly convex functions,
are the ones which cause difficulties (and these are the functions
used to build counterexamples examples below).  Strictly convex functions are generally
easier to approximate, because for these functions, the convexity
constraint is not active.  So the challenge is to build
approximations to convex functions which are robust on nonstrictly convex functions.

\subsection{Failure of naive approaches}
Naive approaches to characterizing discrete convex functions can be
shown to fail.  The first (which naturally corresponds to a finite
element approach) is to fix a triangulation and work with convex
functions on the given triangulation.  This approach fails because it
severely restricts the admissible convex functions to a subset of all
convex functions.  See~\cite{ChoneLeM} for details. 

The second approach is to discretize the condition that the Hessian of a
convex function is positive definite.  This approach, which is more naturally suited to finite
differences,  fails because the discrete Hessian  test can give
false positives and false negatives, see the second example below. 
\begin{example}[Testing convexity in coordinate directions is
  insufficient]
  It is well known that convexity in coordinate directions is not
  enough to ensure convexity.  For example, the function $u(x,y) = xy$
  is linear (and thus convex) in each coordinate direction, but not
  jointly convex in $x$ and $y$.  In particular, the function restricted to the
  direction $y = -x$ is concave.
\end{example}

This last example shows that convexity must be enforced all directions.  The approach of~\cite{CLRM} was to enforce convexity in
all directions available on the grid.  This required a superlinear
(approximately $\mathcal{O}(N^{1.8})$) number of constraints in the
number of grid points $N$.  The resulting constraints were complicated
to implement and computationally intractable.

\begin{example}[The discrete Hessian test fails]\label{exDHF}
  Two examples in~\cite{Morin1} show that the discrete Hessian test
  can fail.  The first is an example of a non-convex grid function,
  which nevertheless has a positive definite discrete Hessian.  The
  second example is a piecewise linear function $u(x,y) = \max( x, ax
  + by, y - c)$ which is convex but its discrete Hessian which is
  not positive definite.
\end{example}

On the other hand, for smooth ($C^3$), strictly convex functions, it
is true that (for small enough $\dx$) the discrete Hessian is positive
definite~\cite[Theorem 4.2]{Morin1}.  The theoretical drawback of this
approach is that polyhedral functions are excluded from the
class of admissible convex functions: instead smoothed and convexified
versions are included.  For example, in their computation of the
Monopolist problem (see \S\ref{monopolist}), the result is a smoothed approximation to the
solution.  Our approach has the advantage that piecewise linear
functions on the grid are admissible and the numerically exact
solution is obtained for the Monopolist example (see~\S\ref{sec:mono}).  The practical
drawback of their approach is that it involves a conic optimization
problem, to enforce the condition that the Hessian is positive
definite.  Our  optimization problems has linear
constraint, which is generally simpler than conic constraints~\cite{BoydBook}.

\begin{example}[{Interpolations of convex functions may not be convex}]
  Given a convex function, the piecewise-linear triangulation of the
  function may not be convex.  The following example comes from~\cite{CLRM}. Consider $u(x,y) = (x+y)^2$ on $[0,1]^2$.  Then the
  piecewise linear approximation of $u$ using the two triangles with
  common edge given by the line $y=x$ is \emph{concave}.
\end{example}

It is also shown in~\cite{CLRM} that, for fixed triangulations,
enforcing convexity condition only on adjacent triangles results in an
admissible set which is strictly smaller than $\C$.

\section{Polyhedral approximations to the cone of convex functions}

In this section, we describe the polyhedral approximation to the cone
of convex functions.  In fact, we bracket the cone of convex functions by a family of interior cones, and a family of exterior cones, which we call  \emph{uniformly convex} and  \emph{mildly nonconvex}, respectively.  In the limit as the directional resolution goes to zero, both of these cones should approach the cone of convex functions.

The method we present uses finite difference discretizations
of variational problems, along a local
characterization of the convexity constraints.  It has the flavor of
PDE methods. Indeed, we use results and estimates first obtained in
solving a PDE for the convex envelope~\cite{ObermanConvexEnvelope,
  ObermanCENumerics}.

Our approach is to approximate the (non-polyhedral) cone of convex
functions by a polyhedral cone which can be represented by linear
inequalities.  This approach leads to simple linear 
constraints: directional second derivative constraints for the continuous problem, and sparse linear constraints 
for the discretized problem.

We give a local characterization of the cone of convex functions.
This results in a constrained variational problem where the number of
constraints is linear in the number of variables.

\subsection{Directional convexity}
The set, $\C$, of convex functions is nonpolyhedral, which is difficult to
enforce numerically.  We approximate this set by the polyhedral
set of directionally convex functions, over a collection of direction
vectors.  Define the directional resolution of a collection of $d$-dimensional unit
vectors $\{ v_i \}, i=1,\dots, k$ as
\begin{align}
\label{dtheta}
\dth &\equiv \max_{\|w\| = 1} \min_i  \cos^{-1}(w^T v_i).
\end{align}
Note that $\dth$ measures the maximum angle between an arbitrary
vector, and the direction vectors for the stencil.  In the two-dimensional case, $\dth$ is half the maximum angle between direction
vectors.  In practise, accurate results are obtained with coarse directional resolution, see Figure~\ref{fig:schemes}.

We can quantify the degree of non-convexity (or uniform convexity)
allowed by the scheme with directional resolution $\dth$.  In the theorem below, we assume for simplicity,  that the function is twice-differentiable. However, the results can be extended to continuous convex functions  by mollification and application of Alexandroff's theorem, as was done in~\cite{EkelandMB}.  Another approach is to use viscosity solutions, as was done in~\cite{ObermanConvexEnvelope}.   

\begin{prop}\label{prop:approxConvex}
Let $u(x)$ be a twice-continuously differentiable function defined on $\R^n$.
Let $V = \{v_i\}_{i=1}^k$ be a set of direction vectors, with directional resolution $\dth$. Assume $\dth \le \pi/4$.  If  
\bq\label{Prop1}
\frac{d^2 u }{dv_i^2} \geq 0, \quad i=1,\dots, k,
\eq
then $u$ is nearly convex, in the sense that
\bq\label{QuantConvexity}
\frac{ \lm}{\lp} \geq -\tan^2(\dth),
\eq
where $\lm \le \dots \le \lp$ are the eigenvalues of the Hessian of $u$ at $x$.
If in addition
\bq\label{Prop2}
\min_{v_i \in V} \left\{ \frac{d^2 u }{dv_i^2} \right\}
\geq  \tan^2(\dth) 
\max_{v_i \in V} \left\{  \frac{d^2 u }{dv_i^2} \right\} 
\eq
then
\bq\label{uconvex}
\lambda_1 \ge 0
\eq
and $u$ is convex.
\end{prop}

\begin{remark}
The conditions~\eqref{Prop1} are linear inequality constraints.  By introducing additional variables, the condition~\eqref{Prop2} can also be implemented as linear constraint.  This is achieved by the linear inequality constraints
\[
0 \le \gamma \le \frac{d^2 u }{dv_i^2} \le \Lambda
\]
and
\[
\gamma \ge \tan^2(\dth) \Lambda.
\]
Strictness in the inequalities can easily be forced, for example by adding the term $\Lambda - \gamma$ to the objective function.

\end{remark}

\begin{proof} Let $u(x)$ be the given twice-differentiable function, and suppose the minimum of $\lm[u](x)$ occurs at $x$.   
Let $w_1$ be the eigenvector corresponding to $\lm$.
Let $\theta$ be the (positive) angle between $w_1$ and the nearest grid direction $v_i$.  By~\eqref{dtheta},  
$\theta \leq \dth$.
Decompose $v_i = \cos\theta w_1 + \sin\theta w$, where $w$ is a unit vector orthogonal to $w_1$.    Then compute
\begin{align*}
\frac{d^2 u }{dv_i^2} &= (\cos\theta w_1 + \sin\theta w)^T D^2u\, (\cos\theta w_1 + \sin\theta w)\\
& = \cos^2\theta w_1^T D^2u\, w_1 + \sin^2\theta w^T D^2u\, w + 2 \sin\theta\cos\theta w_1^T D^2u\, w
\\
& = \cos^2\theta \lm + \sin^2\theta w^T D^2u\, w 
\\
&\le 
  \cos^2\theta \lm + \sin^2\theta \lp. 
\end{align*}
In the computation above, we have used the fact that  $w_1$ is an eigenvalue  and $w$ is orthogonal to $w_1$ to eliminate the cross term, and we have also used the estimate $ w^T D^2u\, w  \le \lp$.

Now if \eqref{Prop1} holds, the previous calculation yields 
\[
\frac{\lm}{\lp}  \ge  -\tan^2(\dth),
\]
since $\theta \leq \dth$, and~\eqref{QuantConvexity} is established.

Now, performing a similar calculation to the one above, but setting $w_n$ to be the eigenvalue corresponding to $\lp$, and $v_k$ to be the nearest grid direction to $w_n$, we obtain
\begin{align*}
\frac{d^2 u }{dv_k^2} 
&\ge   \cos^2\theta \lp + \sin^2\theta \lm. 
\end{align*}

Assuming~\eqref{Prop2} holds, we can combine the previous two inequalities with \eqref{Prop2} to obtain
\begin{align*}
\cos^2\theta \lm + \sin^2\theta \lp
&  \ge \tan^2(\dth)  \left (\cos^2\theta \lp + \sin^2\theta \lm \right)
\\
&\ge \tan^2(\theta)  \left (\cos^2\theta \lp + \sin^2\theta \lm \right)
\end{align*}
using the fact that $\tan^2(\dth) \ge \tan^2(\theta)$.  Simplifying gives
\[
\cos^4 \theta \lm \ge \sin^4\theta \lm
\] 
Recall the assumption that  $\theta < \dth < \pi/4$, which means that $\cos(\theta) \ge \sin(\theta)$.  Thus, $\lm \ge 0$ and~\eqref{uconvex} is established.
\end{proof}

\section{Presentation of the variational problem}

Let $\Omega$ be a closed bounded convex subset of $\R^n$, and
\[
\C \equiv \{  u : \Omega \to \R \mid u \text{ is convex in } \Omega \}.
\]
Consider the variational problem subject to the convexity constraint
\bq
\inf_{u \in \C \cap \K }  
I(u), \quad\text{where}\quad I(u) \equiv \int_\Omega f(x,u(x), \grad u(x))\, dx,
\eq and $\K$ is a closed convex subset of a given space $X =
H^1(\Omega)$ or $X = L^2(\Omega)$.  Here $\K$  represents, for example, boundary conditions or pointwise bounds. 

For the purposes of the numerical optimization, we require that for
each fixed $x$, $f$ is either a quadratic function of $u$ and $\grad
u$, or $f$ is convex and homogeneous of order one.  In the first case,
after expressing the approximation of $\C$ by linear constraints, we
arrive at a quadratic program (QP), and in the second case we arrive at
linear program (LP).  Both are standard convex optimization problems, and can
be solved by commercial or academic software packages~\cite{BoydBook}.

If $I$ is lower semicontinuous, coercive, and strictly convex on $\K$,
the existence and uniqueness of a minimizer directly follows from
standard arguments.

\subsection{Example problems}

\begin{example} A typical example is given by
\[
I_f(u) = \int_\Omega \frac 1 2 | \grad u  |^2 + fu \,d\x
\]
subject to $u \in \C$, and possibly with Dirichlet boundary conditions.
\end{example}

\begin{example}
  The Rochet-Chon\'e example,~\cite{RochetChone},  is given by
\[
I_{RC}(u) = \int_\Omega \frac 1 2 (u_x^2 + u_y^2) + x u_x + y u_y - u \, dx \, dy
\]
subject to $u \ge 0$, and $u \in \C$, with no explicit boundary conditions.
\end{example}

\begin{example}
A standard problem is the projection (in some norm) onto the set of
convex functions, which is given by the minimizer of the functional
\[
I_{u_0}(u) = \| u - u_0\|_{X},
\]
where $X = H^1(\Omega)$ or $X = L^2(\Omega)$, subject to $u \in \C$.
\end{example}

\begin{example} The problem of finding convex envelopes of a given
  function is also included (see~\cite{CLRM}).   Given $u_0 \in
  L^2(\Omega)$, take $X = L^2(\Omega)$ and
\[
I_{u_0}(u) =  \| u - u_0\|_{L^2(\Omega)}
\]
along with the additional constraint
\[
\{ u\in L^2(\Omega) \mid u \le u_0 \text{ a.e.} \}.
\]
The convex envelope can also be computed directly using a
PDE~\cite{ObermanConvexEnvelope,ObermanCENumerics}.
\end{example}

\subsection{An analytic solution in one dimension}
In this section we give an explicit non-trivial analytic solution in
one dimension.  The analytic solution is useful for verification of the numerical method.  In this one-dimensional case, the solution is the convex envelope
of the minimizer of the unconstrained problem, as was shown in~\cite[\S
4]{CLRM}.  However, as was also shown in~\cite[\S
4]{CLRM}, in higher dimensions this is no longer the case.

Consider 
\[
\min J[u] = \int_{-1}^1 \frac 1 2 u_x^2(x) + f(x) u(x) dx
\]
over functions in $H^1$ with Dirichlet boundary conditions $u(-1) =
u(1) = 0$.  The convexity constraint becomes simply $ u_{xx} \ge 0$.
Here we restrict to the special case of a step function, 
\[
f(x) = 
\begin{cases}
-c & x \leq 0\\
+c & x \geq 0
\end{cases}
\]
for some $c > 0$.
The minimizer will solve $u_{xx}(x) = f(x)$ on some unknown set, and will solve $u_{xx} = 0$ outside.   Furthermore, it's easy to see that the set on which $u_{xx} > 0$ will be of the form $[a, 1]$, for some $a \ge -1$.  Thus the candidate solutions will be convex, differentiable functions which are 
 linear on  $[-1,a]$ and quadratic on~$[a,1]$.

\subsubsection*{Description of the admissible set}
By imposing the differentiability restriction at  $x = a$, 
we can write candidates for the minimizer as a one-parameter family
\begin{equation}\label{ua}
u^a(x) = 
\begin{cases}
m(x+1) & \text{ for }  x\in [-1,a]\\
c(x-b)(x-1)/2 & \text{ for } x \in [a, 1]
\end{cases}
\end{equation}
where 
$m = -{c}(a-1)^2/4$,  and $b =  (a^2 + 2a -1)/2$.
Note that each $u^a(x)$ is a local minimizer of the functional $J$.  

\subsubsection*{Evaluation of $J[u]$}
Now compute $J[u^a]$ as a function of $a$.  
The result, arrived using elementary integration, is:
\[
J[u^a(x)] = - \frac 1 {48} c^2 (a - 1)^2  (3 a^2  + 10 a - 1).
\]
Now we can find a local minimum of $J$, by setting $\frac{d}{da} J[u^a(x)] = 0$, to arrive at
\[
\frac 1 4 c^2 (1 -a) (a^2  + 2 a - 1) = 0
\]
which gives the minimum at $
a = \sqrt{2} -1.
$

\subsubsection*{Verification of the Euler-Lagrange equation}
The solution of the unconstrained problem is
\[
u(x) =
\begin{cases}
-f(x)x(x+1)/2 & \text{ for } x < 0\\
+f(x)x(x-1)/2 & \text{ for } x \geq 0.
\end{cases}
\]
The convex hull of the solution is obtained by putting $b=0$ into \eqref{ua}, which gives $a^2  + 2 a - 1=0$ and agrees with the solution.

\begin{figure}
\scalebox{.4}{
\includegraphics{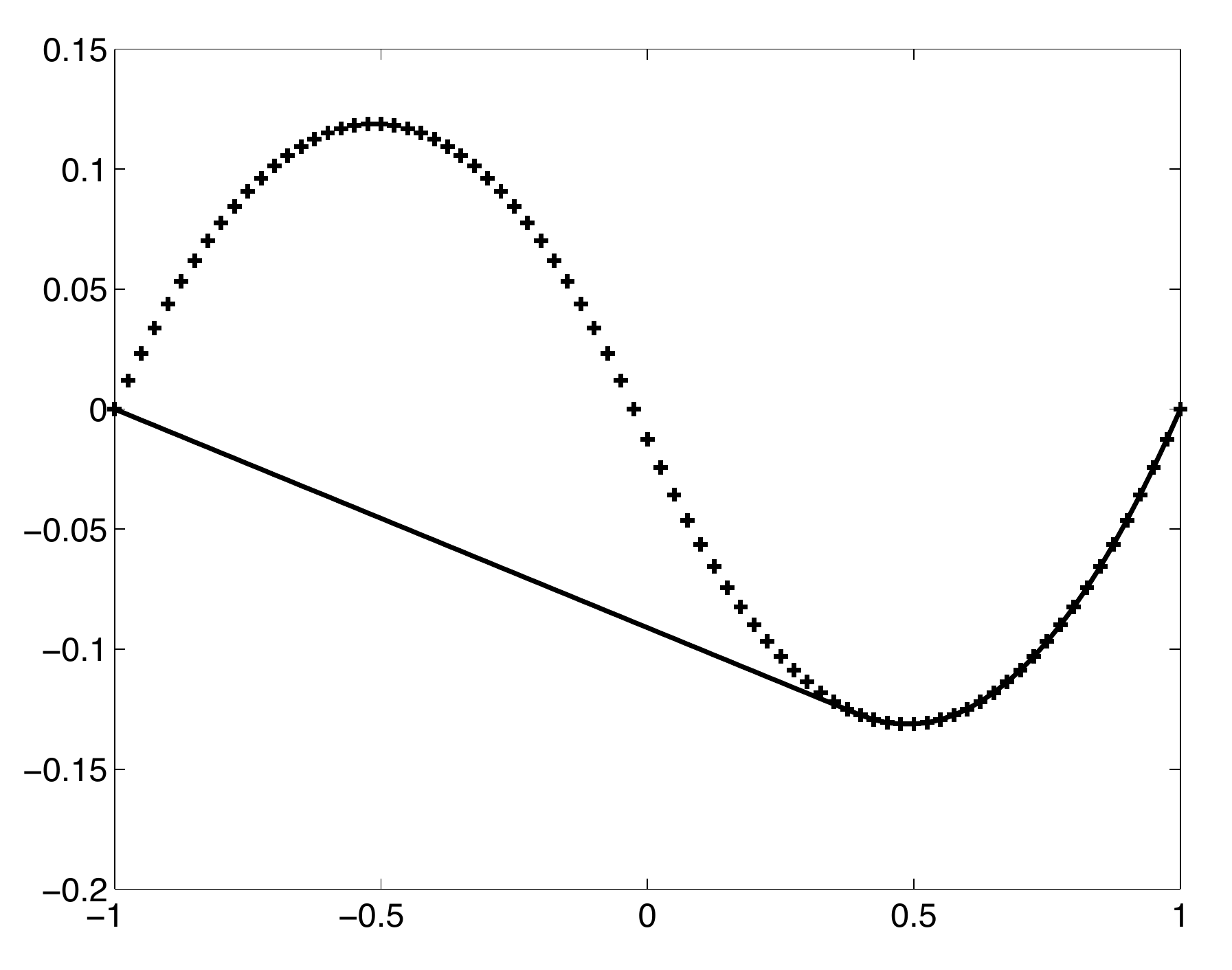}
}
\caption{
Minimizer of the constrained problem (solid line), and the minimizer of the unconstrained problem.  
}
\label{fig1}
\end{figure}

The global minimizer, along with the solution of the unconstrained
problem is plotted in Figure~\ref{fig1}.

\subsection{The Monopolist problem}\label{monopolist}
The Monopolist problem of~\cite{RochetChone} is given by the
variational problem
\[
\min_{u\in C} I[u] = \int_\Omega c \abs{\grad u}^2/2 - u(\x) + \grad
u(\x)\cdot \x \, d\x
\]
for $\Omega =  [0,1]^2 \subset \R^2$
subject to 
\[
u \ge 0.
\]
The exact solution is given in~\cite{RochetChone}.  The quantity of interest in the economics problem is the gradient map, which determines the optimal production choice of the monopolist as a function of the distribution of the population, which depends on two parameters.  In this case, the gradient map  is a combination of a constant, linear functions, and nonlinear functions.  These parts of the map correspond to a production of a fixed good for a large part of the population, and varying degrees of customization for the other parts of the population.

\subsection{A variant of the Monopolist problem}\label{monopolistvar}
A variant of the Monopolist problem was studied in~\cite{MMdesign} and 
computed by~\cite{Morin1}.  It has the special feature of a linear objective function, which was easier to compute using the methods of~\cite{Morin1}. The problem is given by
\[
\min_{u\in C} I[u]  =  \int_\Omega  u(\x) -  \grad u(\x)\cdot \x \, d\x
\]
for $\Omega =  [0,1]^2 \subset \R^2$
subject to 
\[
0 \le u_x, u_y \le 1, \quad  u(0) = 0.
\]
The analytical solution of this problem is given by
\[
u(x,y) = \max( 0, x - a, y -a, x + y - b ), \quad 
a = \frac 2 3, \quad b = \frac 1 3 (4 - \sqrt 2),
\]
with optimal value 
\[
I(u) = \frac 2 {27} ( 6 + \sqrt 2 ) \approx .0549.
\]

\section{Approximation and implementation}
\subsection{Overview of the  discretization of the problem}
Consider the functional
\[
I^n(u) = \sum_i f(x_i,u(x_i), \grad u(x_i))\, dx,
\]
where $\grad u(x_i)$ is a finite difference approximation.  The
discretization of the gradients in the objective function is performed
using standard centered finite differences.  Significant performance
differences resulted from different discretizations. 

The convexity constraints are implemented via a set of linear
inequalities, at each point, which are directional second derivatives:
\[
u(x_i) \leq \frac{ u(x_i + h v_j) + u(x_i - h v_j) }{2}, 
\quad (i = 1,\dots, n),~(j = 1,\dots, k)
\]
where $v_j$ is a collection of \emph{direction vectors}.

Recall that $\K$ is a closed convex subset of a given space $X =
H^1(\Omega)$ or $X = L^2(\Omega)$, and represents, for example, boundary conditions or pointwise bounds.  We assume that $\K$  can be represented by linear inequalities.

The result is a quadratic (or linear) minimization, with linear constraints.

\subsection{Convergence of convex approximations}
Proofs of convergence of approximations of variational problems with convexity constraints can be found in~\cite{CLRM} and in~\cite{EkelandMB}.  
In this article, we focus on the structure of the resulting discrete optimization problems produced by the approximation to the cone of convex functions. 

\section{Numerical Discretization}
In order to have a well-posed (and accurate) convex variational
problem, we need to carefully build the constraint matrices.  We found
that specific choices of these constraint and objective function matrices
make a significant difference in the solution time for the variational
problem.  The choice of boundary conditions is also important.
For simplicity, we present small  examples of the matrices used to generate the computational results.

\subsection{One dimensional constraint matrices}
In one dimension, convexity is enforced only at interior nodes.
The convexity constraint  operator $D_{xx}$ is given for $n=5$ by:
\[
D_{xx} = \frac{1}{\dx^2} 
\left[\begin{array}{rrrrr}1 & -2 & 1 & 0 & 0 \\0 & 1 & -2 & 1 & 0 \\0 & 0 & 1 & -2 & 1\end{array}\right].
\]

\subsection{One dimensional objective function gradient matrices}
It is important that the operator corresponding to $|u_x|^2$ in the objective function is a square symmetric matrix which is zero on constant functions.  

This is accomplished  as follows.
Start with the forward and backward finite difference operators, on the largest possible stencil. For example, in one dimension on a small grid,
\[
D^+_{x} = \frac{1}{\dx} 
\left[\begin{array}{rrrrr}
    -1    & 1  &  0 &    0   &  0\\
     0   & -1  &   1   &  0   &  0\\
     0    & 0 &   -1   &  1  &  0 \\
     0    & 0&     0    & -1 &   1\\
\end{array}\right].
\]
The operator corresponding to $u_x^2$ in the objective function is given by
\[
 \frac{1}{\dx}
\left (
D^+_x + (D^+_{x})^T
\right )
=  (D^+_{x})^T D^+_{x}
= \frac{1}{\dx^2} 
\left[\begin{array}{rrrrr}
    1    &-1    & 0 &    0    & 0  \\
    -1    & 2  &  -1 &    0   &  0\\
     0   & -1  &   2   & -1   &  0\\
     0    & 0 &   -1   &  2  &  -1\\
     0    & 0&     0    &-1 &    1\\
\end{array}\right]
\]

\subsection{Two dimensional objective function gradient matrices}
The two dimensional versions of the operator corresponding to $|\grad u|^2$ can be constructed as a square symmetric matrix in a similar way.

The gradient squared objective function corresponds to a symmetric, positive definite Laplacian operator.  The right way to build this operator with finite difference matrices is to average the full forward and backward difference operators in each grid direction.  This ensures that the final operator is symmetric.  (In fact corner values can be included by doing this with diagonal operators if this is desired).

So we build
\[
-\left (D_{xx} + D_{yy} \right) 
= \frac{1}{\dx}
\left \{
D^+_x + D^-_x + D^+_y + D^-_y
\right \},
\]
where the operators correspond to the one dimensional case above in
each coordinate.  Thus, in the $n=3$ case, the symmetric positive
Laplacian matrix is the $9\times 9$ matrix
\[
-D_{xx + yy} = \frac{1}{h^2}
\left [ \begin{array}{rrrrrrrrr}
     2  &  -1 &    0   & -1  &   0&     0  &   0&     0  &   0\\
    -1  &   3 &   -1   &  0  &  -1 &    0  &   0 &    0 &    0\\
     0  &  -1 &    2   &  0  &   0 &   -1  &   0   &  0 &    0\\
    -1  &   0  &   0 &    3  &  -1 &    0  &  -1    & 0  &   0\\
     0  &  -1   &  0&    -1   &  4 &   -1&     0&    -1 &    0\\
     0   &  0  &  -1  &   0  &  -1  &   3    & 0&     0 &   -1\\
     0  &   0   &  0   & -1  &   0   &  0  &   2  &  -1 &    0\\
     0   &  0 &    0    & 0 &   -1   &  0 &   -1  &   3 &   -1\\
     0    & 0&     0    & 0&     0  &  -1&     0   & -1&     2\\
 \end{array}\right ].
\]

\subsection{Two dimensional constraint function gradient matrices}
If the gradient matrix appears in the constraints, we can afford to use higher accuracy centered differences away from the boundary, and we use 
 forward or backward differences near the
boundary.  This which yields
\[
D_{x} = \frac{1}{2\dx} 
\left[\begin{array}{rrrrr}
    -2    & 2    & 0 &    0    & 0  \\
    -1    & 0  &  1 &    0   &  0\\
     0   & -1  &   0   & 1   &  0\\
     0    & 0 &   -1   &  0  &  1\\
     0    & 0&     0    &-2 &    2\\
\end{array}\right].
\]

\section{Accuracy of the convexity constraint}
Proposition~\ref{prop:approxConvex} gives a characterization of the approximation to the cone of convex functions.  In this section we perform numerical tests to measure the accuracy of this approximation.

\begin{figure}
\centering
\subfigure[Grid directions for the width 1, width 2 schemes.]{\scalebox{.3}{\includegraphics{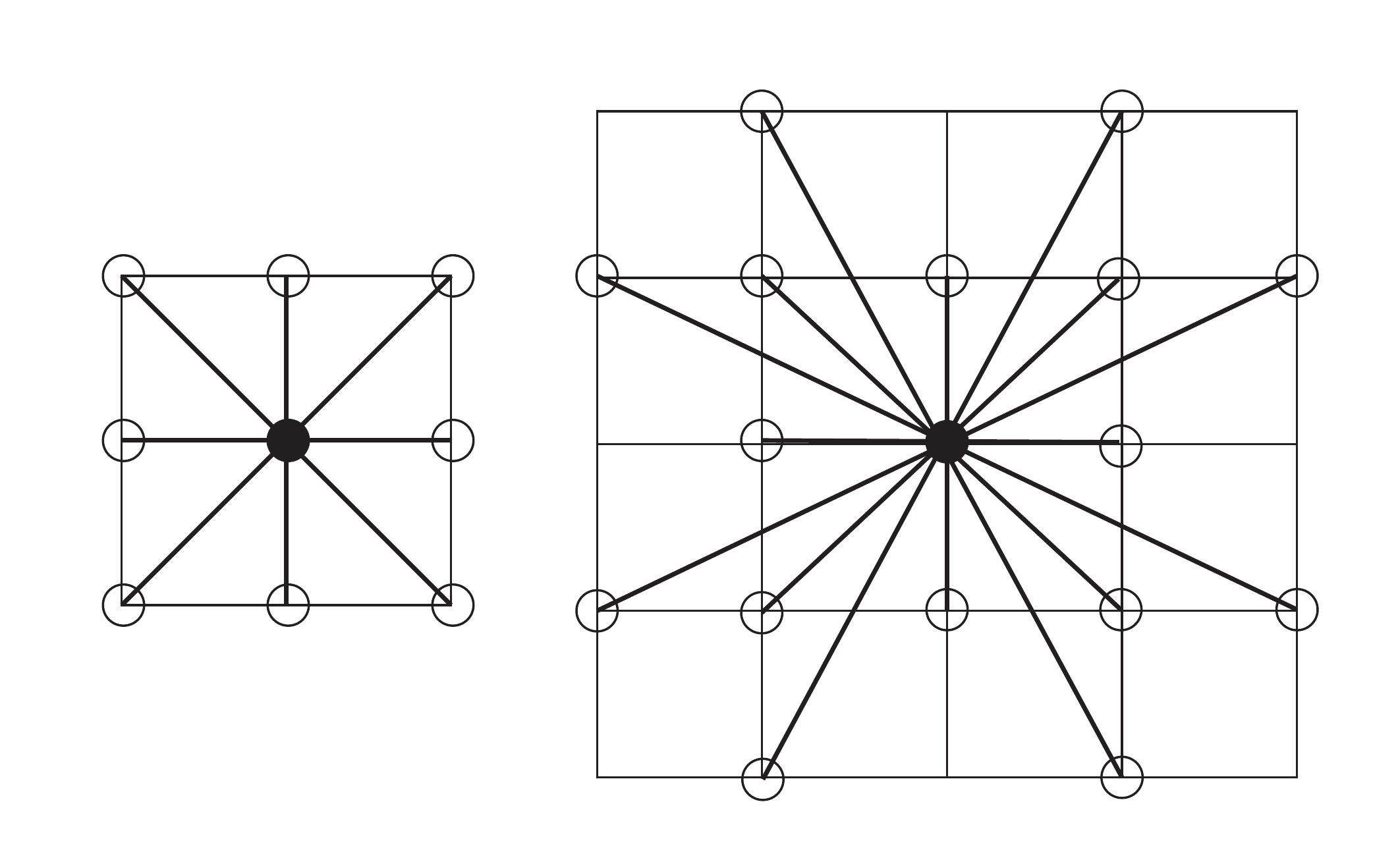}}}
\subfigure[Grid directions for the width 1, 2, 3, 4 schemes.]{\scalebox{.3}{\includegraphics{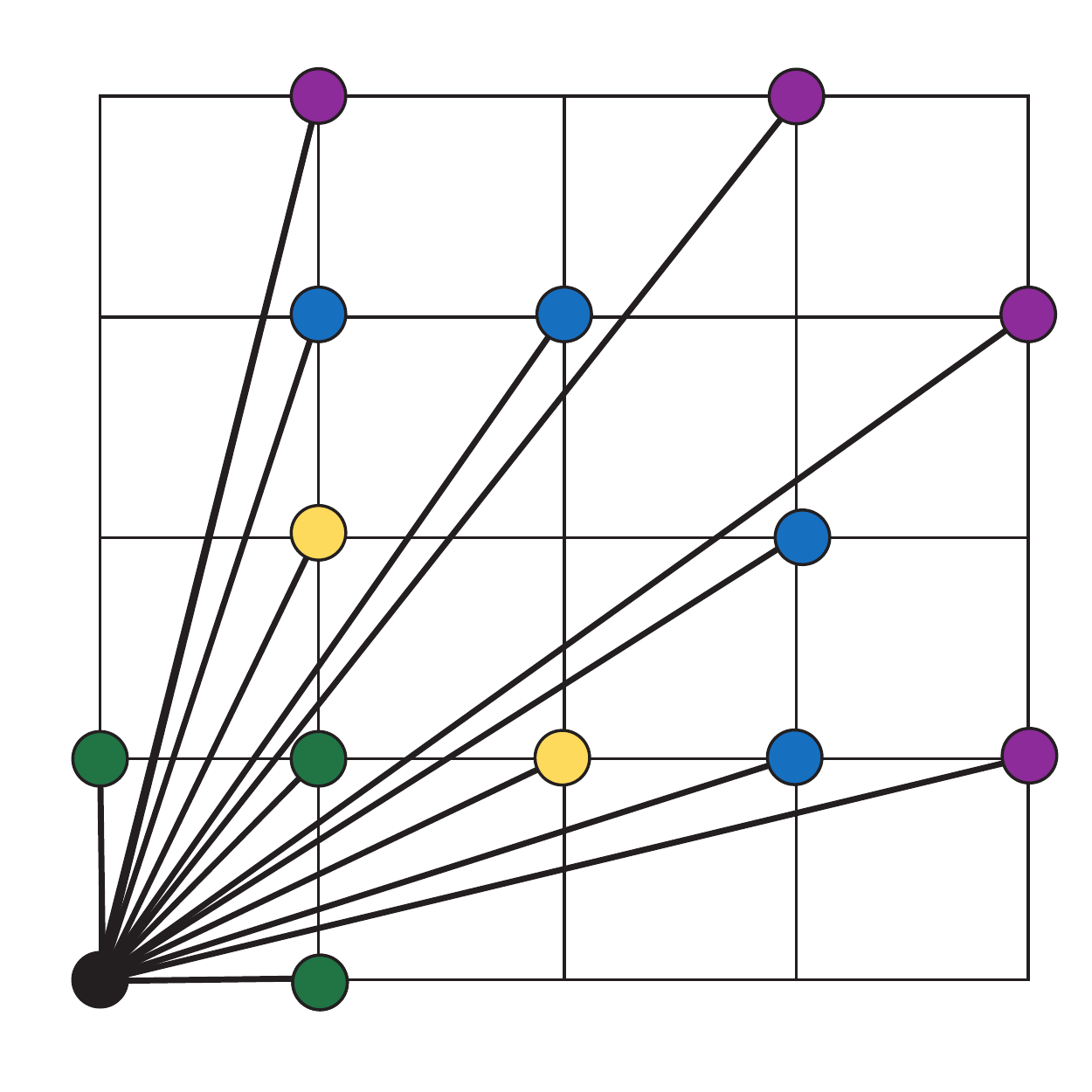}}}
\caption{Illustration of the directions in which convexity is tested for different schemes.  Only directions in the first quadrant are shown in (b). }\label{fig:schemes}
\end{figure}

\begin{table}
\begin{center}
\begin{tabular}{cllllll}
Width $(w)$ &  \multicolumn{4}{l}{(Additional) direction vectors} & $\dth = \tan^{-1}(1/w)/2$ & $\tan^2(\dth)$ \\ 
\hline
1		& (1,0) 	& (1,1) 	& 		&		& .39  & .17\\
2        	& (2,1) 	&(1,2) 	&		&		& .23  & .056\\
3         	& (1,3) 	& (3,1)     	& (2,3) 	&(3,2)   	& .16  &  .026\\
4               & (1,4)		& (4,1)		& (3,4)   	&(4,3)       & .12  &  .015         
\end{tabular}
\caption{Directions vectors in the first quadrant for schemes of width $w$.  Values of $\dth$ and $\tan^2 \dth$.
}\label{table:nbrs}
\end{center}
\end{table}
   
\subsection{Comparison of Wider Stencil Convexity Constraints}
In this section we investigate the effect of enforcing convexity in more directions.
For a fixed grid size of  $21\times 21$, we compared the effect of enforcing convexity in the four directions (horizontal, vertical, and the two diagonals) on the 9 point nearest neighbor grid, and the additional 4 directions available on the wider grid.  These grid directions are presented in see Figure~\ref{fig:schemes} and Table~\ref{table:nbrs}.  For the purpose of comparison, the $L^2$ projection was performed on various convex and non-convex test functions.  

For more examples, the worst case predictions of Proposition~\ref{prop:approxConvex} were overly pessimistic.  Indeed, for many functions we tested, the minimizer was the same up to numerical accuracy for the different methods.   For the following functions, the minimizer was the same, for all width $1,2,3$ schemes:
\[
\sin(2\pi x),\quad -x^2, \quad (x-3y)^2, \quad \abs{x-3y}.
\]
Notice that the schemes are invariant on convex functions, but may be (wrongly) invariant on slightly nonconvex functions.

Next we present examples which were sensitive to the number of directions in which convexity was enforced.

\subsubsection*{Results with a spiky nonconvex function}
For other functions, the results with more directions differed strongly.
For example, (always on $[0,1]^2$)
\[
\exp
\left(-30\left(
(x-.5)^2 + (y-.5)^2
\right)
\right)
\]
gave results which converged as the stencil widened.  
The error, measured as the maximum value between the width 4 solution and the width 1,2,3 respectively,
 was $0.02, 0.01, 0.003$, for $n=21$,
 and $0.02, 0.01, 0.005$, for $n=31$.




\subsubsection*{Results with a nonconvex quadratic}
The next example uses the nonconvex quadratic function, 
$(x^2 - \alpha y^2)/2$ rotated by an angle $\theta$:
\begin{multline}
g^{\alpha,\theta}(x,y) = 	\\
	\left (\frac{\cos^2\theta + \alpha \sin^2\theta}{2}\right ) {x^2} 
	+ \left (	(1-\alpha)\cos\theta\sin\theta \right ) xy 
 	+ \left (  \frac{ \alpha\cos^2\theta + \sin^2\theta}{2} \right ) {y^2},
\end{multline}
with $\alpha < 1$, so that $\lm = \alpha$, $\lp = 1$, and the eigenvectors are at an angle of $\theta$ from the coordinate axes.    

The worst case angles, and threshold values of $\alpha$  for detection of non-convexity are given by $\dth/2$, $\alpha$ from Table~\ref{table:nbrs}.
\begin{align*}
\theta &= \pi/8,  && \alpha = -0.055, && \text{ for width 1 } 
\\
\theta &=\pi/12,  && \alpha = -.0264, && \text{ for width 2}
\end{align*}

Taking a simple example, nonconvex quadratics, $x^2 - .5 y^2$ on $[-1,1]^2$, rotated by $0,\pi/2, \pi/4$ radians.  The minimizer in this last case was the same for both schemes of width 1 and 2.

Results were consistent with the predictions of Proposition~\ref{prop:approxConvex}.
However  values of   $\theta = \pi/8$, $\alpha = -.5$, yielded the same minimizers.

For $\theta = \pi/8$, $\alpha = -.1$, the 4 direction scheme failed to see the non-convexity of the function and returned it as a minimizer, while the wider scheme gave something more convex.  The results differed by $.04$.

With $\theta = \pi/8$, $\alpha = -.05$ the solutions of width 1 and width 2 were the same.   For width 3, the solution differed by .018.


%
%
%
 \section{Example Problems and Numerical Results}

\subsection{Projections $d=1$}
In this section, we perform projections in dimension $d=1$.
While the more interesting projections are in two dimensions, this example is useful for performance benchmarking and for visualization.

The function to be projected is
\[
f(x) = \sin(\pi x)
\]
on $[-1,1]$.

\begin{figure}[htdp]
\subfigure[$L^1$ projection]{\includegraphics[width=.49\textwidth]{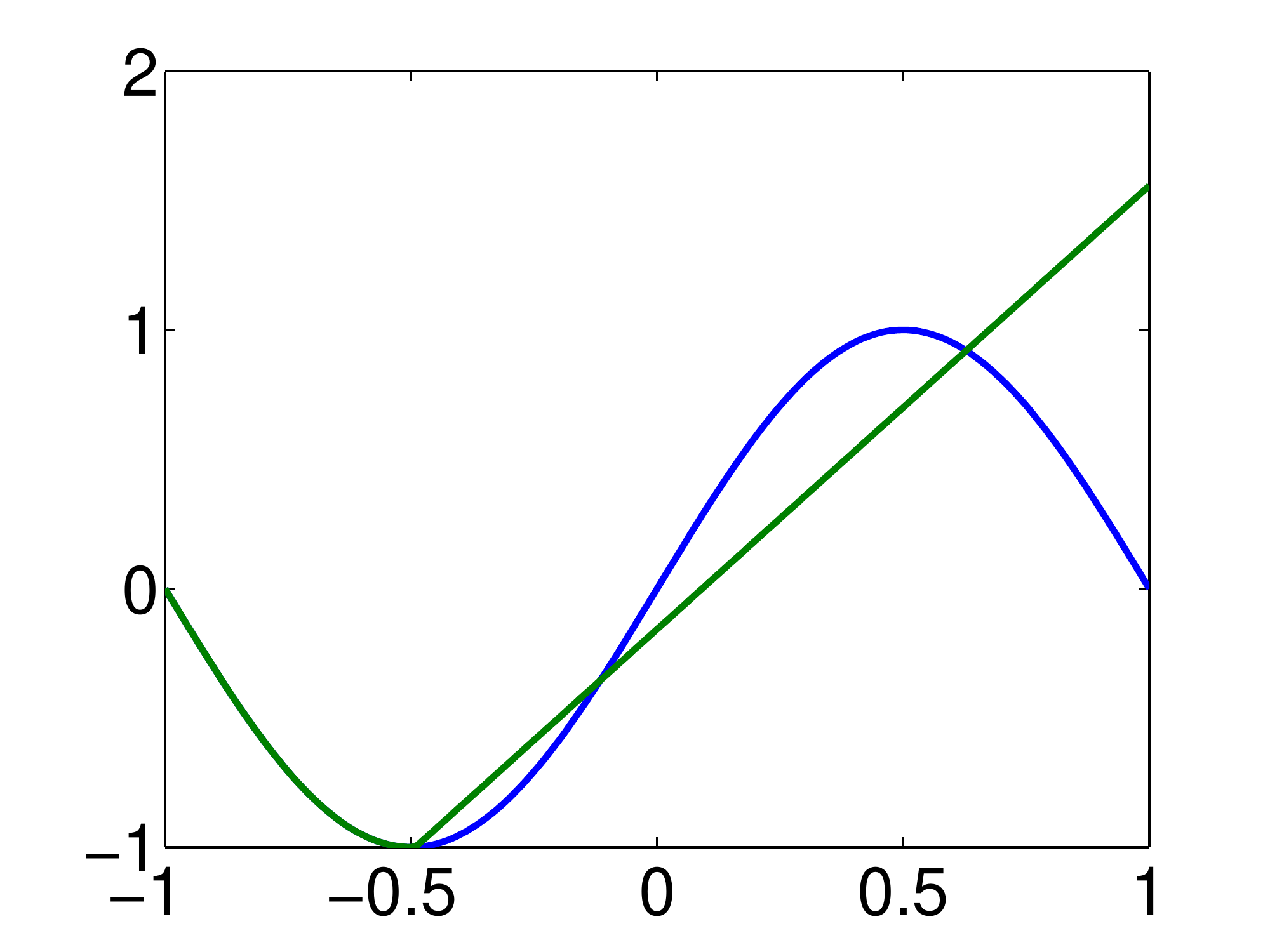}}
\subfigure[$L^2$ projection]{\includegraphics[width=.49\textwidth]{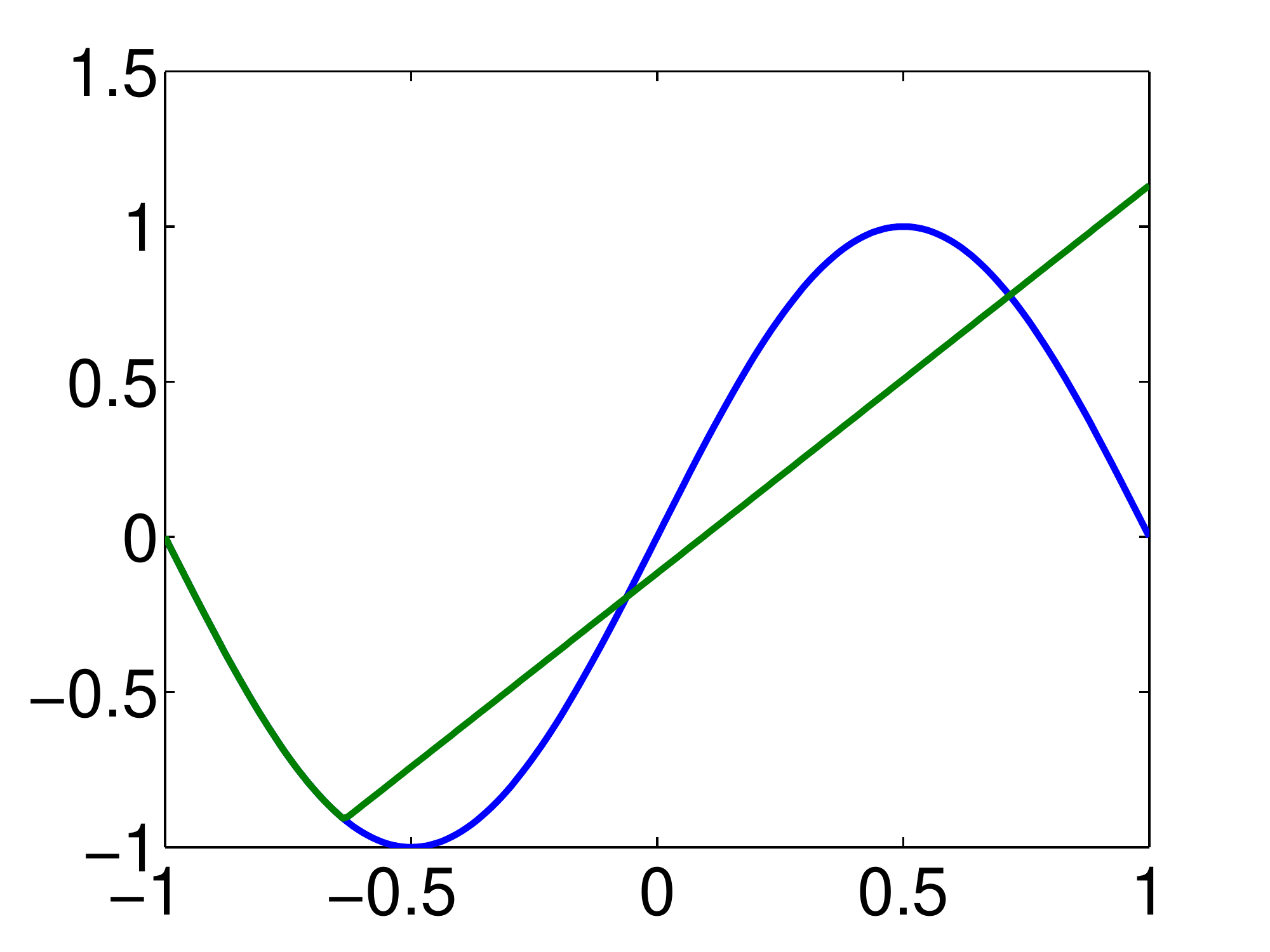}}
\subfigure[$L^\infty$ projection]{\includegraphics[width=.49\textwidth]{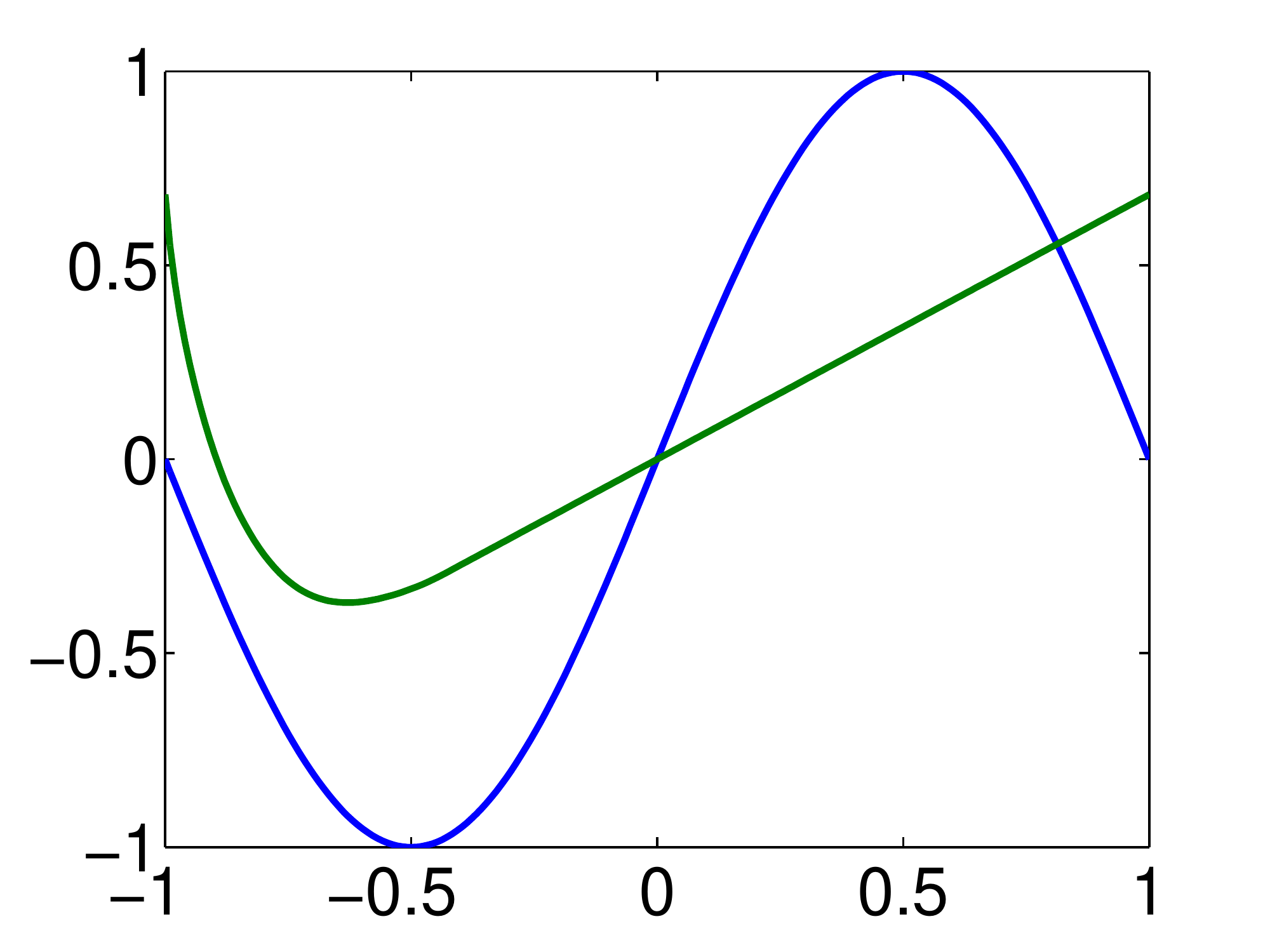}}
\subfigure[$H_0^1$ projection]{\includegraphics[width=.49\textwidth]{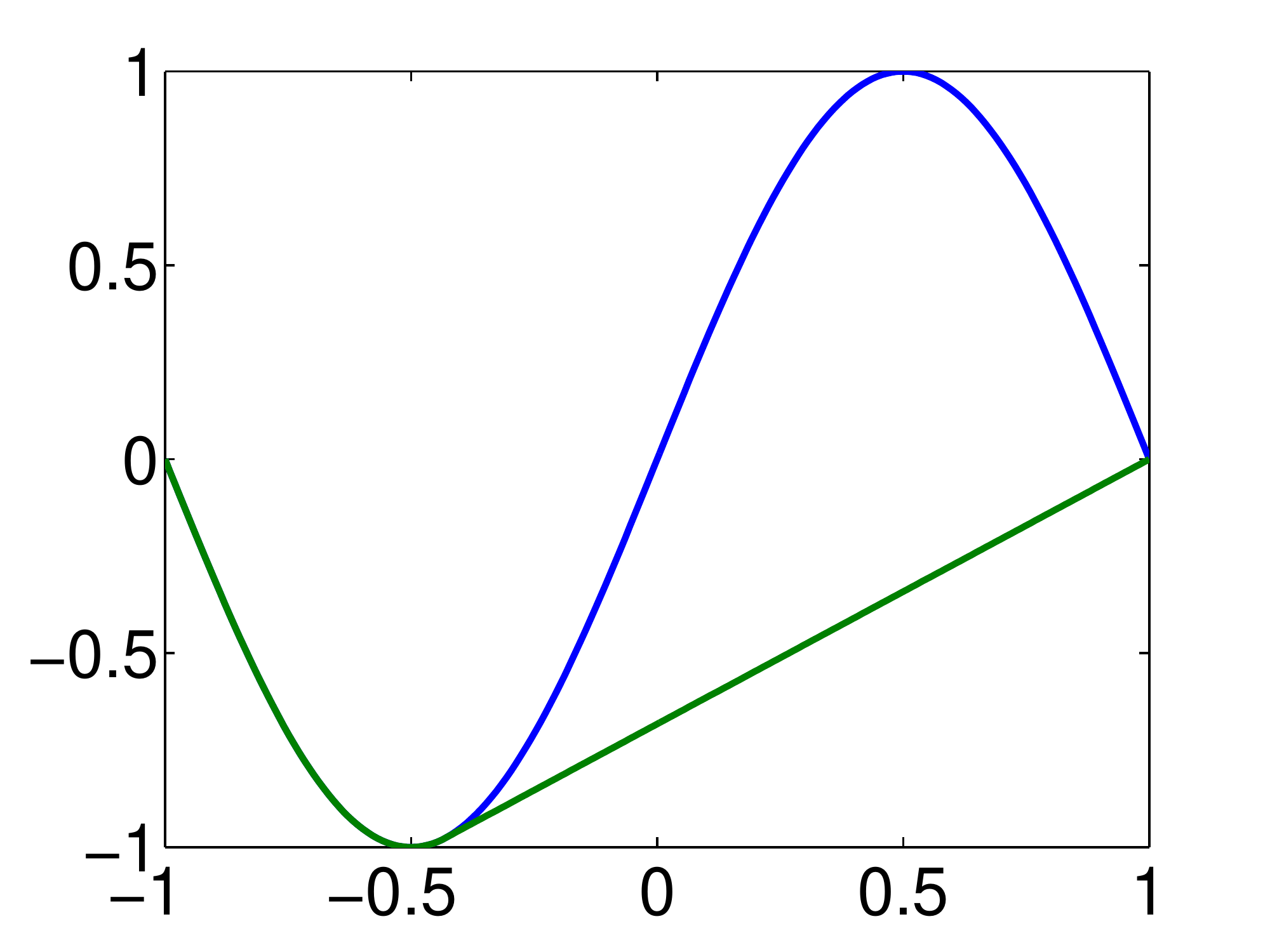}}
\subfigure[$H^1$ projection]{\includegraphics[width=.49\textwidth]{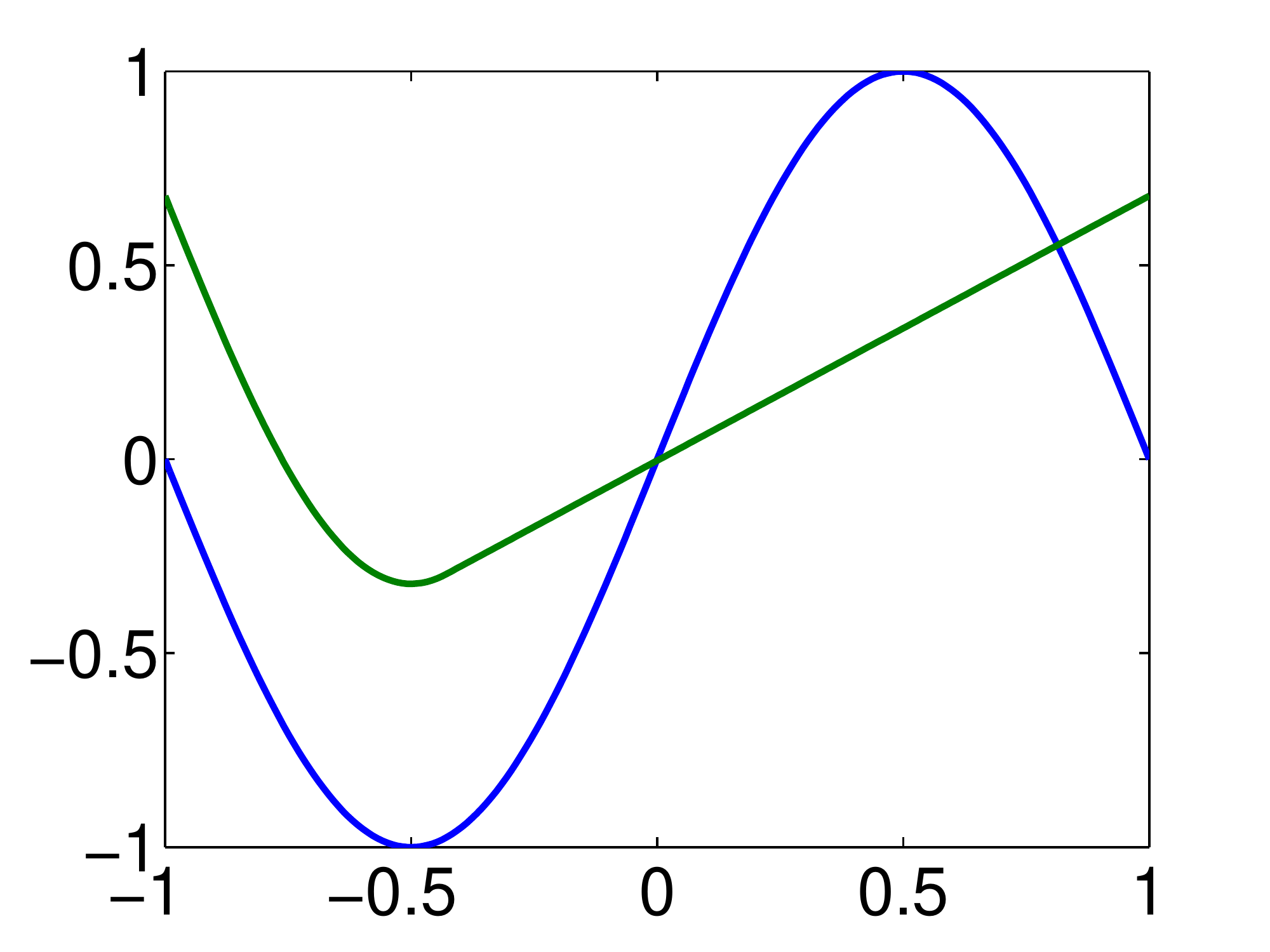}}
\subfigure[$H^1$ with gradient constraints]{\includegraphics[width=.49\textwidth]{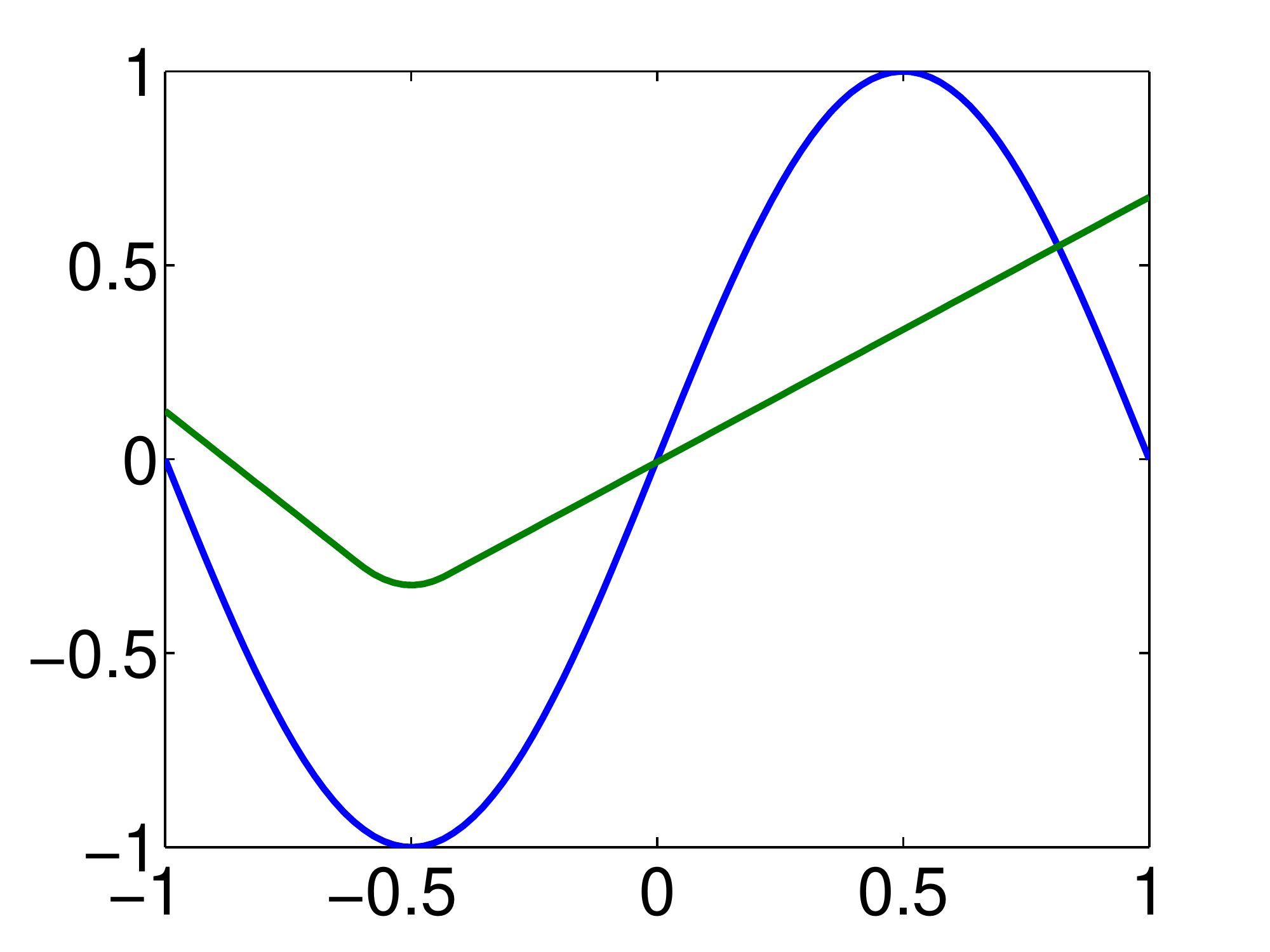}}
\caption{Projection in various norms of $\sin(\pi x)$  onto the cone of convex functions.}
\label{Proj1d}
\end{figure}

Projections in $L^1, L^2, L^\infty, H^1, H^1_0,$  and $H^1$ with gradient constraints of the function $\sin(\pi x)$ are presented in Figure~\ref{Proj1d}.  When needed, the constraint $u(0) = 0$ was also included.

\subsection{Projections $d=2$}
\begin{figure}[htdp]
\subfigure[The function]{\includegraphics[width=.49\textwidth]{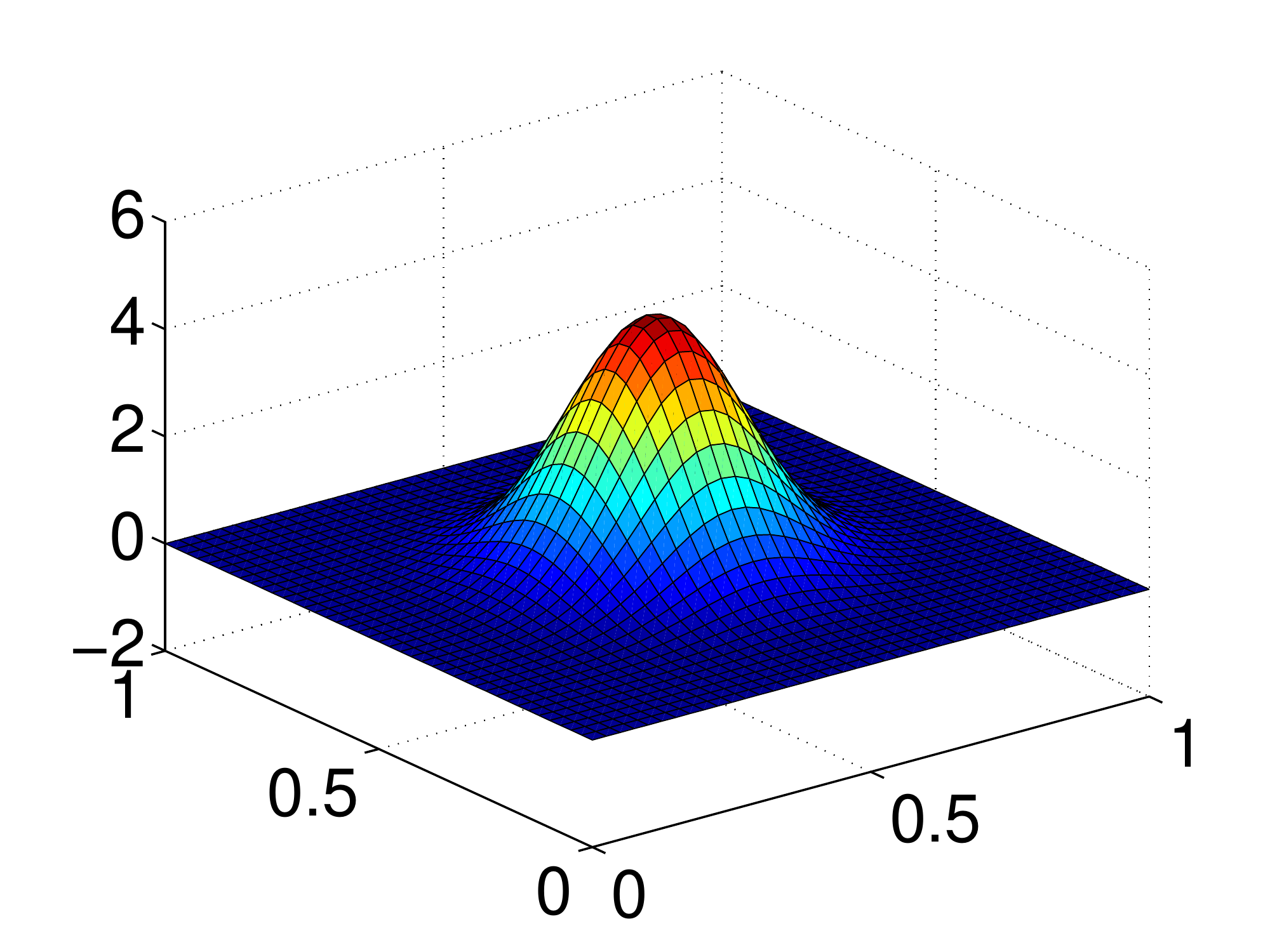}}
\subfigure[$L^1$ projection]{\includegraphics[width=.49\textwidth]{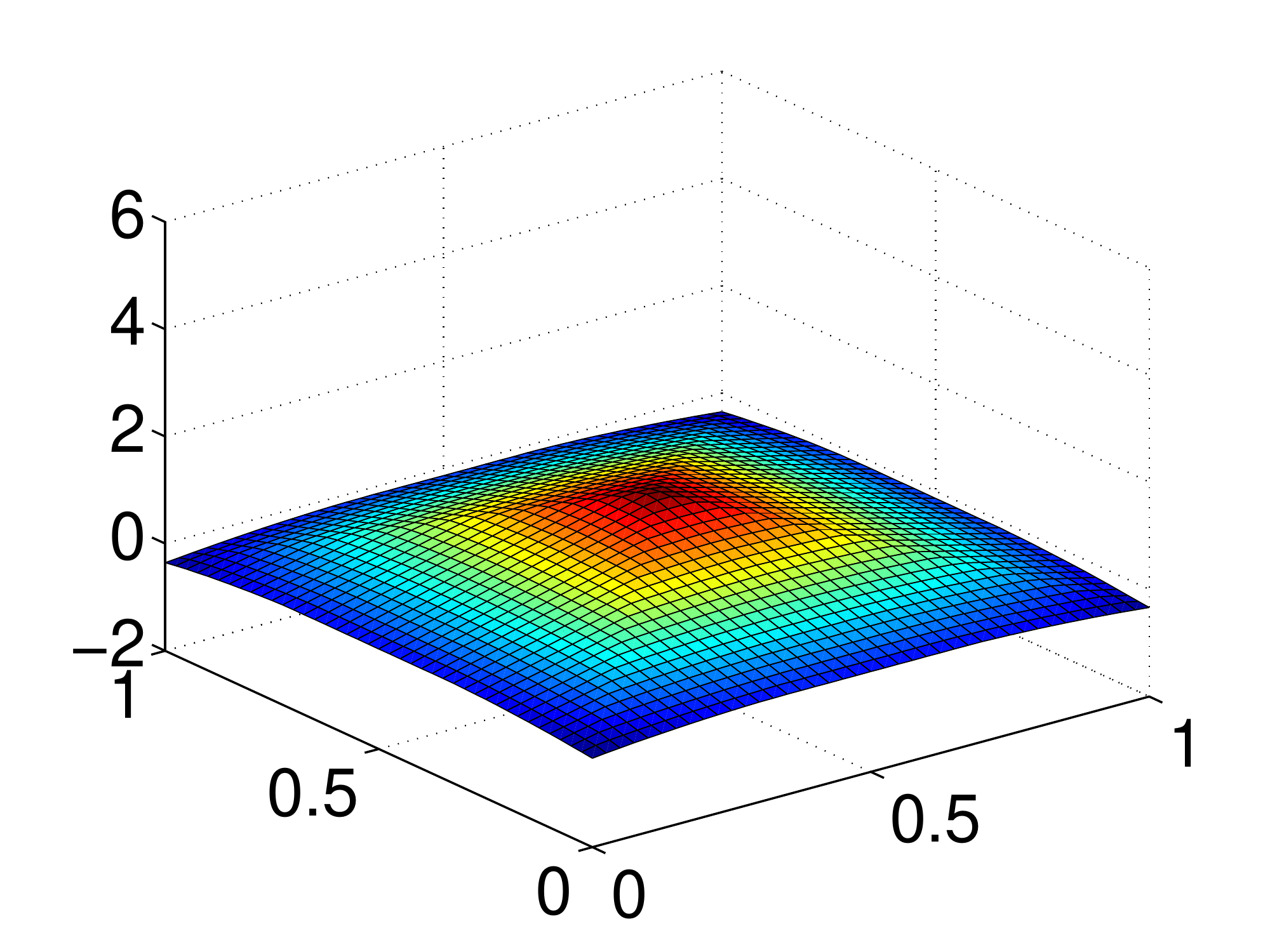}}
\subfigure[$L^2$ projection]{\includegraphics[width=.49\textwidth]{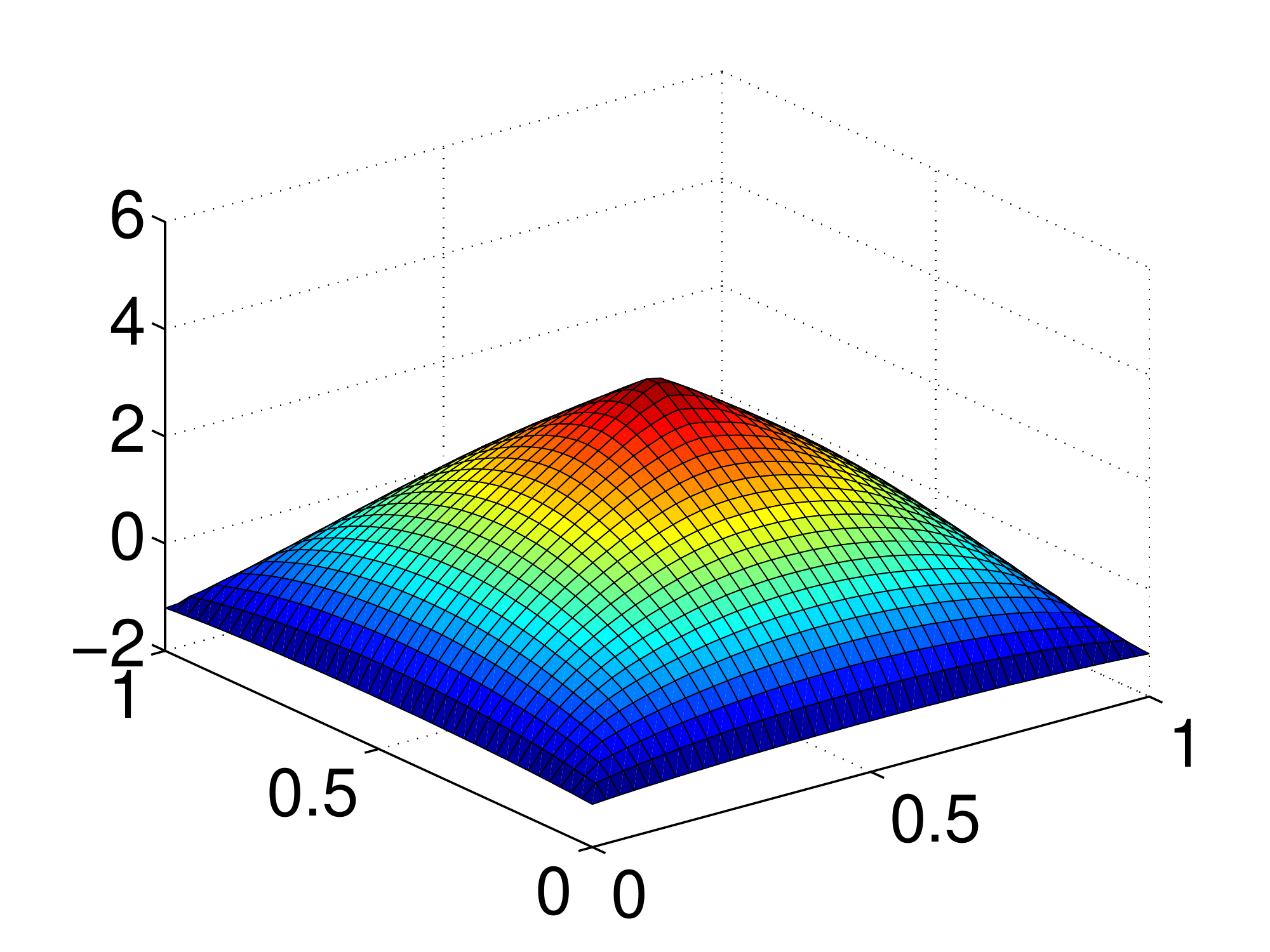}}
\subfigure[$L^\infty$ projection]{\includegraphics[width=.49\textwidth]{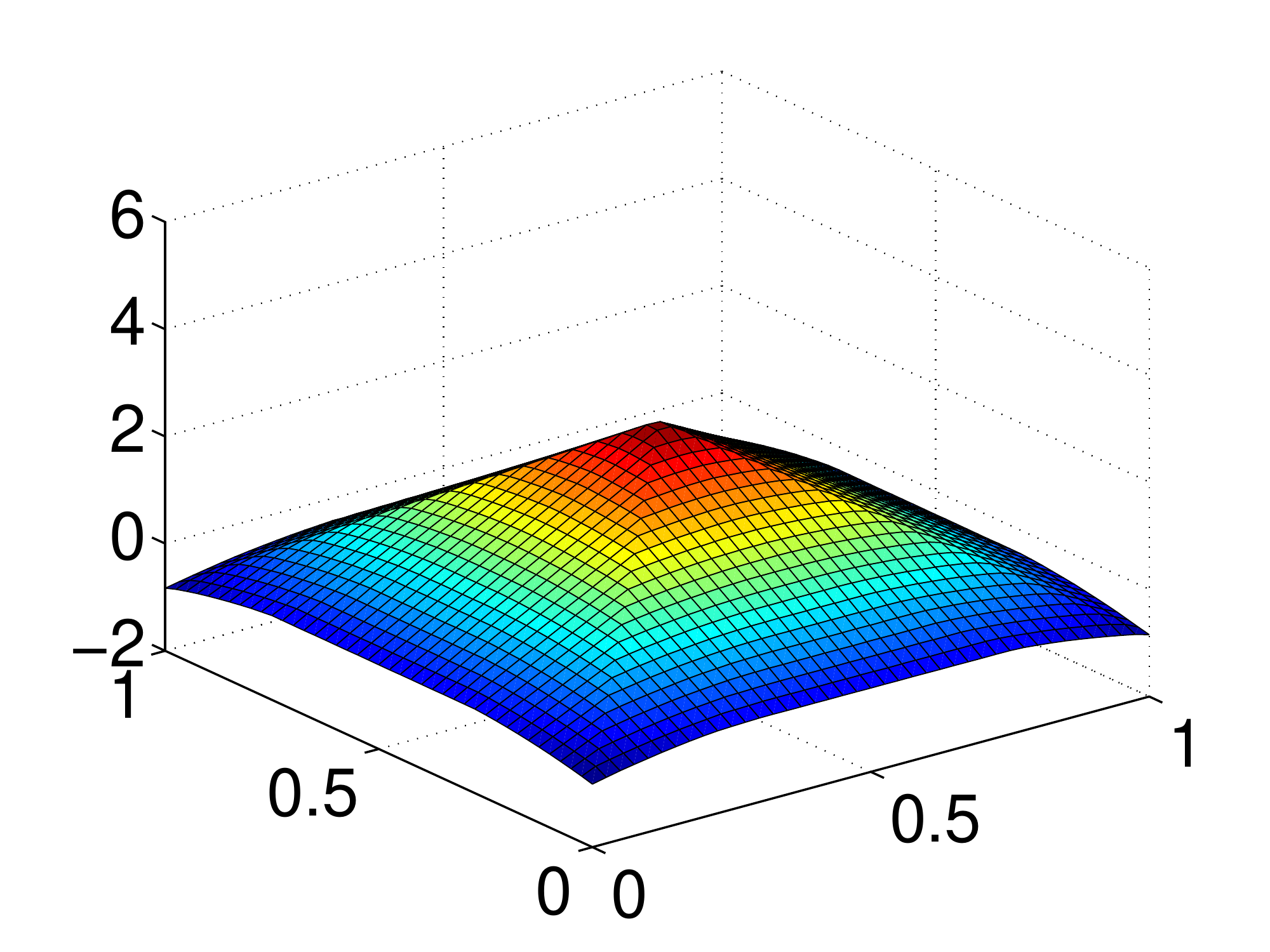}}
\subfigure[$H^1$ projection]{\includegraphics[width=.49\textwidth]{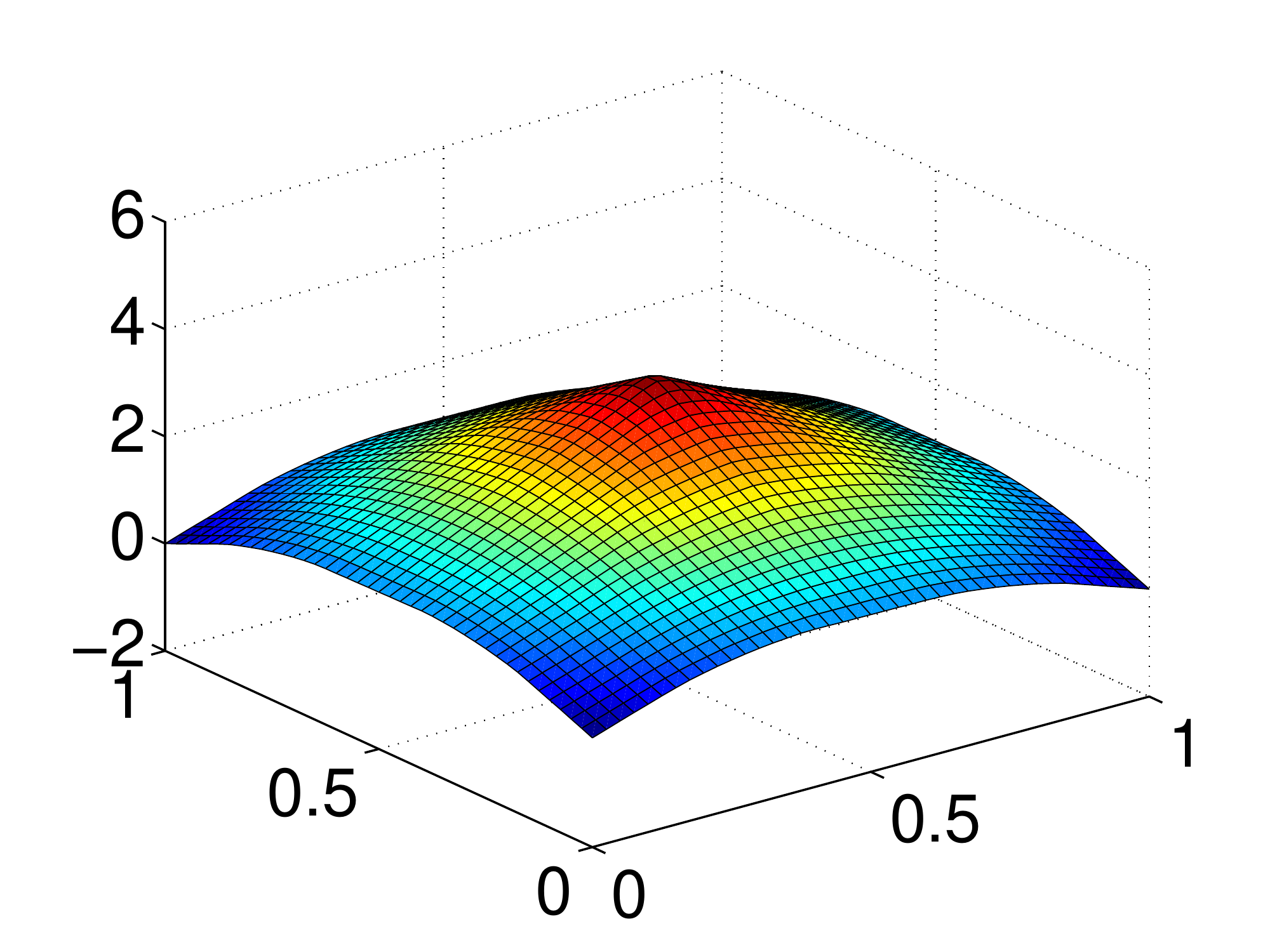}}
\subfigure[$H_0^1$ projection]{\includegraphics[width=.49\textwidth]{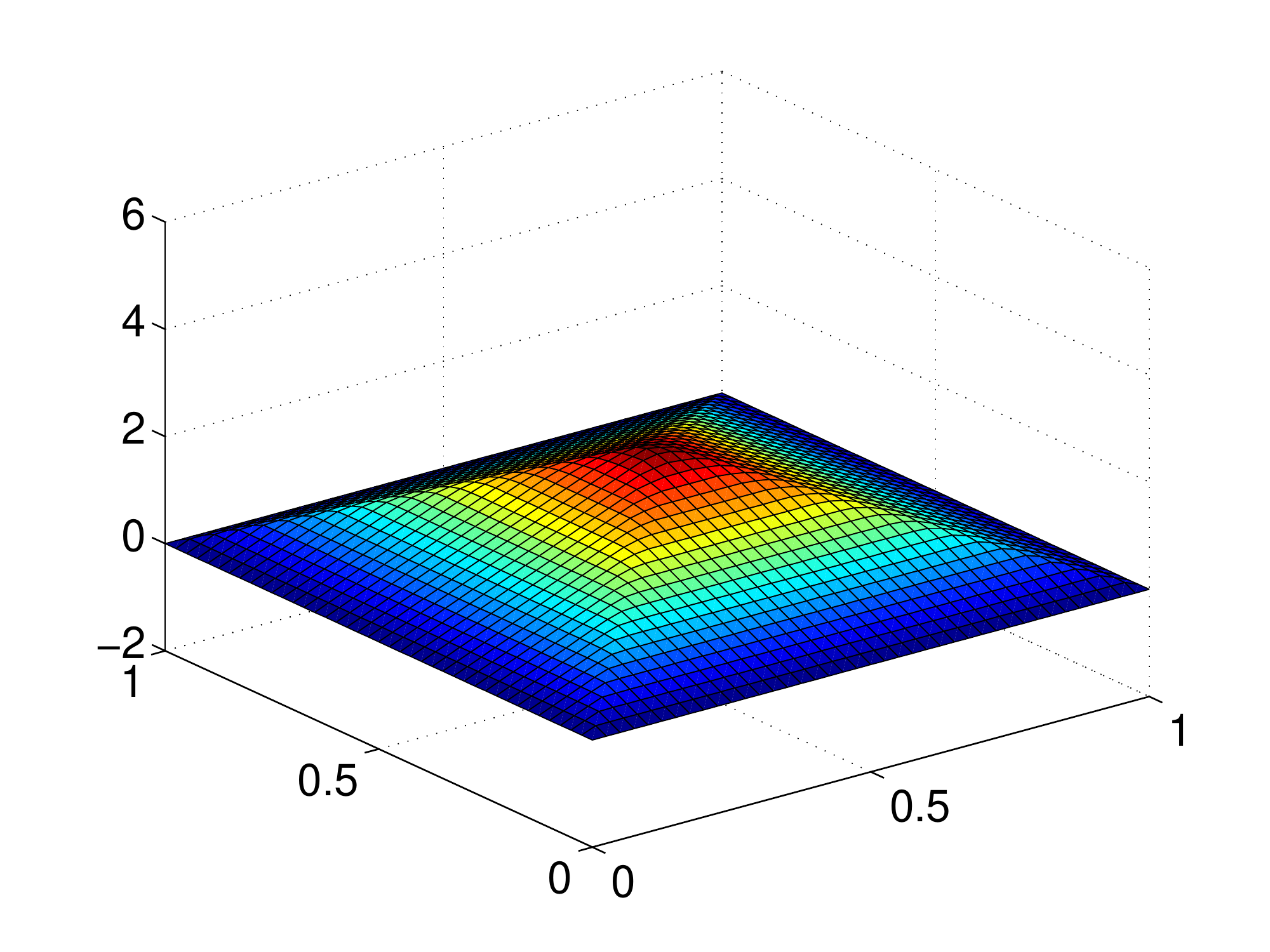}}
\caption{Projection in various norms of the function given by \eqref{spiky} onto the cone of convex functions. }

\label{figProjections}
\end{figure}

In this section, we compute  for comparison purposes  the projections from~\cite{CLRM} and~\cite{Morin1}.  These are projections is various norms with convexity constraints. The function to be projected is given by
\bq\label{spiky}
f(x,y) = -(4 + 5xy^2)\exp
\left(-30\left(
(x-.5)^2 + (y-.5)^2
\right)
\right).
\eq
Plots of the original function and the $L^1, L^2, L^\infty, H^1,$ and $H^1_0$ projections are presented in Figure~\ref{figProjections}.

\subsection{Numerical results for the Monopolist problem}
The example of the Monopolist problem of \S\ref{monopolist} was computed.
The solution, along with a histogram of the gradient map is given in Figure~\ref{fig:mon1}.
The numerical solution captures the fact that the density of the map is concentrated at point, along a line, and then spread out over a quadrant.
\begin{figure}[htdp]
\subfigure[Contour plot of the solution]{ \includegraphics[width=.9\textwidth]{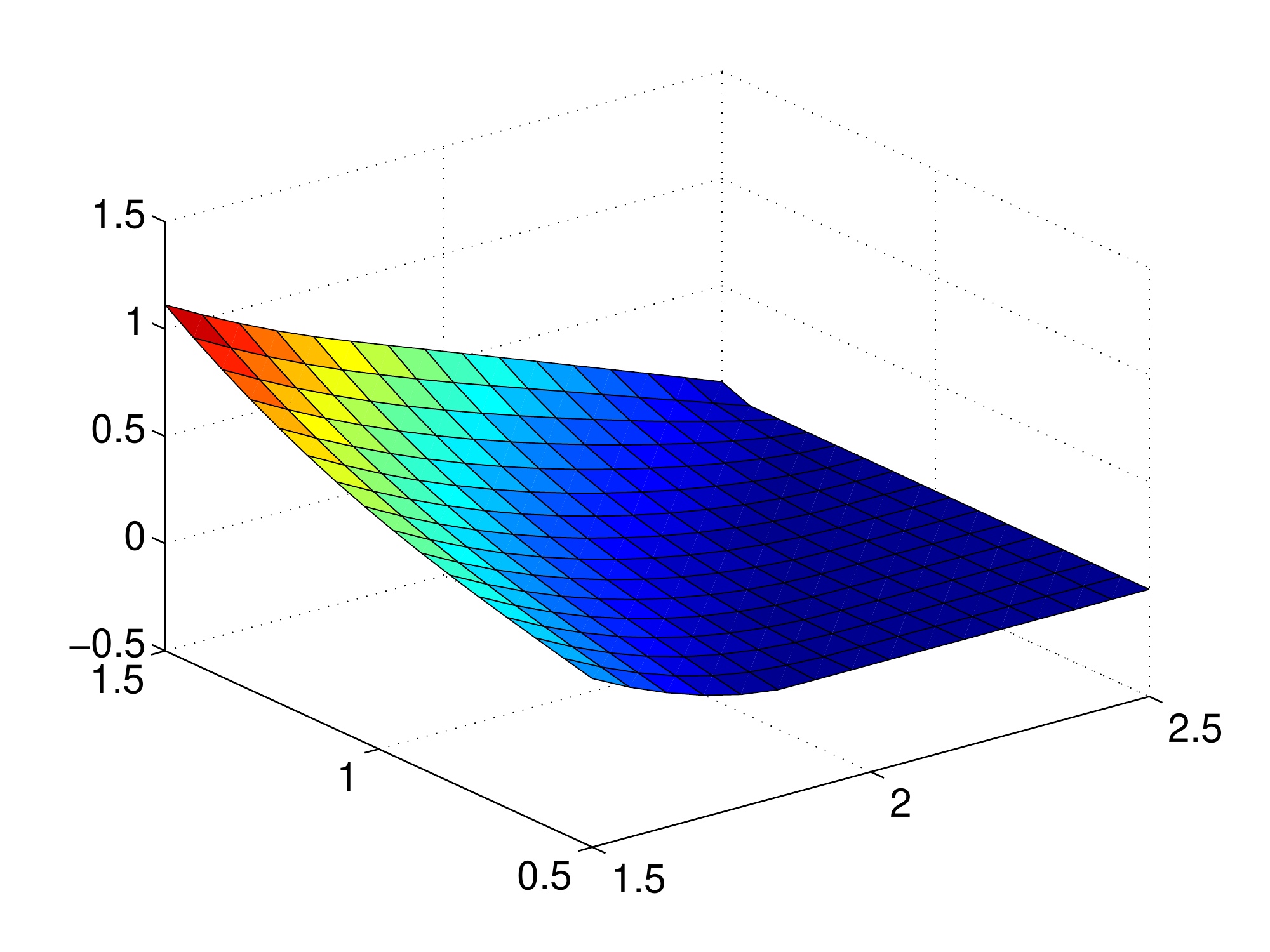}}
\subfigure[Histogram of the gradient]{ \includegraphics[width=.9\textwidth]{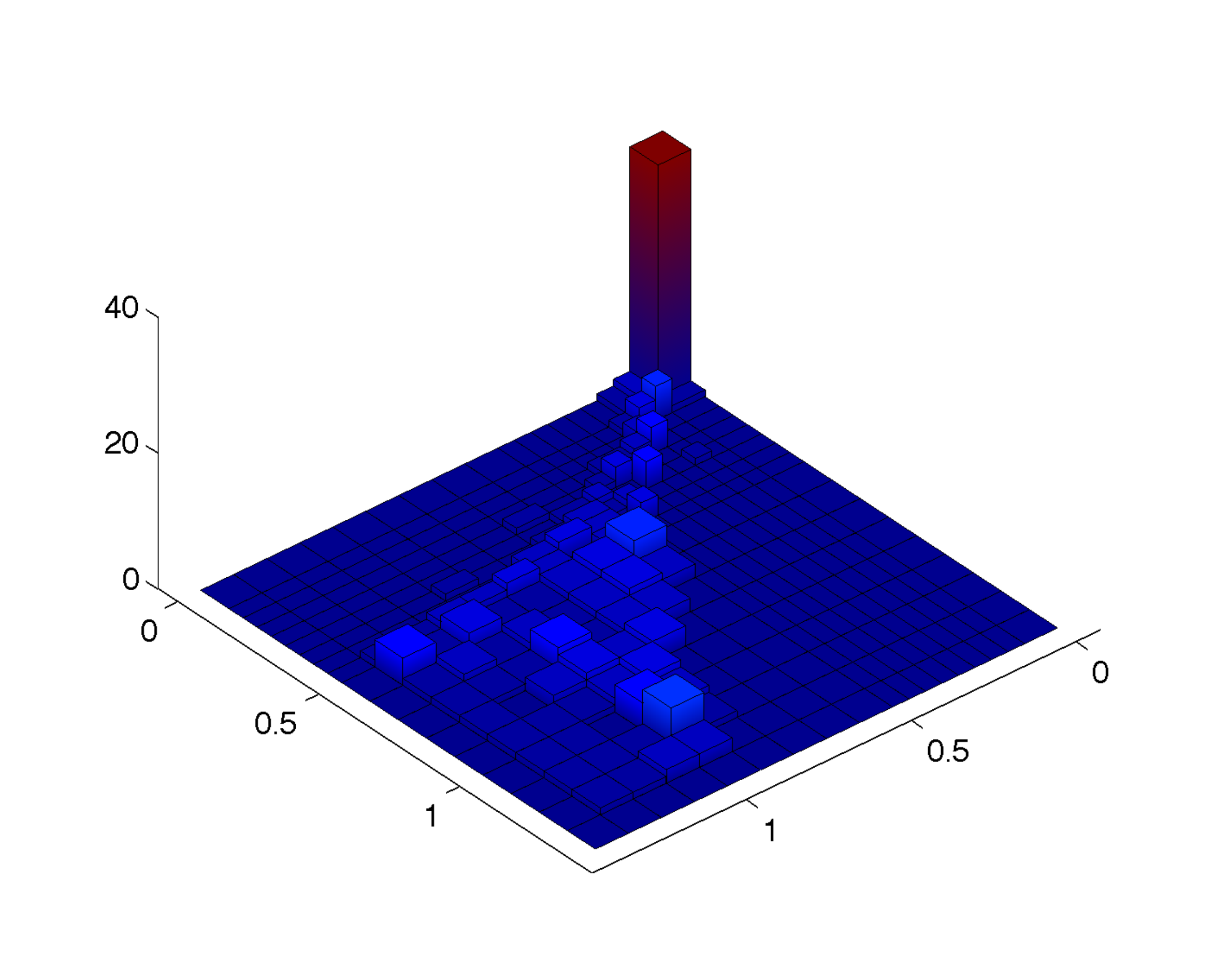}}
\caption{Solution of the Monopolist Problem.
}
\label{fig:mon1}
\end{figure}

\subsection{Numerical results for the variant of the Monopolist problem}\label{sec:mono}
The example of the Monopolist problem of \S\ref{monopolistvar} is well suited to our methods, 
 since our representation of convex functions includes the piecewise linear functions aligned with the coordinate and diagonal directions, which is the case for this example.
 We obtain the exact solution up to interpolation error, and the computational time is which is hundreds of times less than in~\cite{Morin1}.
 
The contour lines of the numerical solutions are presented in Figure~\ref{fig:mon1}.  At this resolution, the numerical solution is indistinguishable from the exact solution.  (This is in contrast to the solution in~\cite{Morin1}, where the contour lines were curved).
The computation was performed using piecewise constant quadrature for the integral, and  the trapezoidal rule.  In both cases only the most narrow (width 1) stencil was used to enforce convexity.
In both cases,  the piecewise linear solution was captured to within interpolation error.

\begin{figure}[htdp]
\includegraphics[width=.7\textwidth]{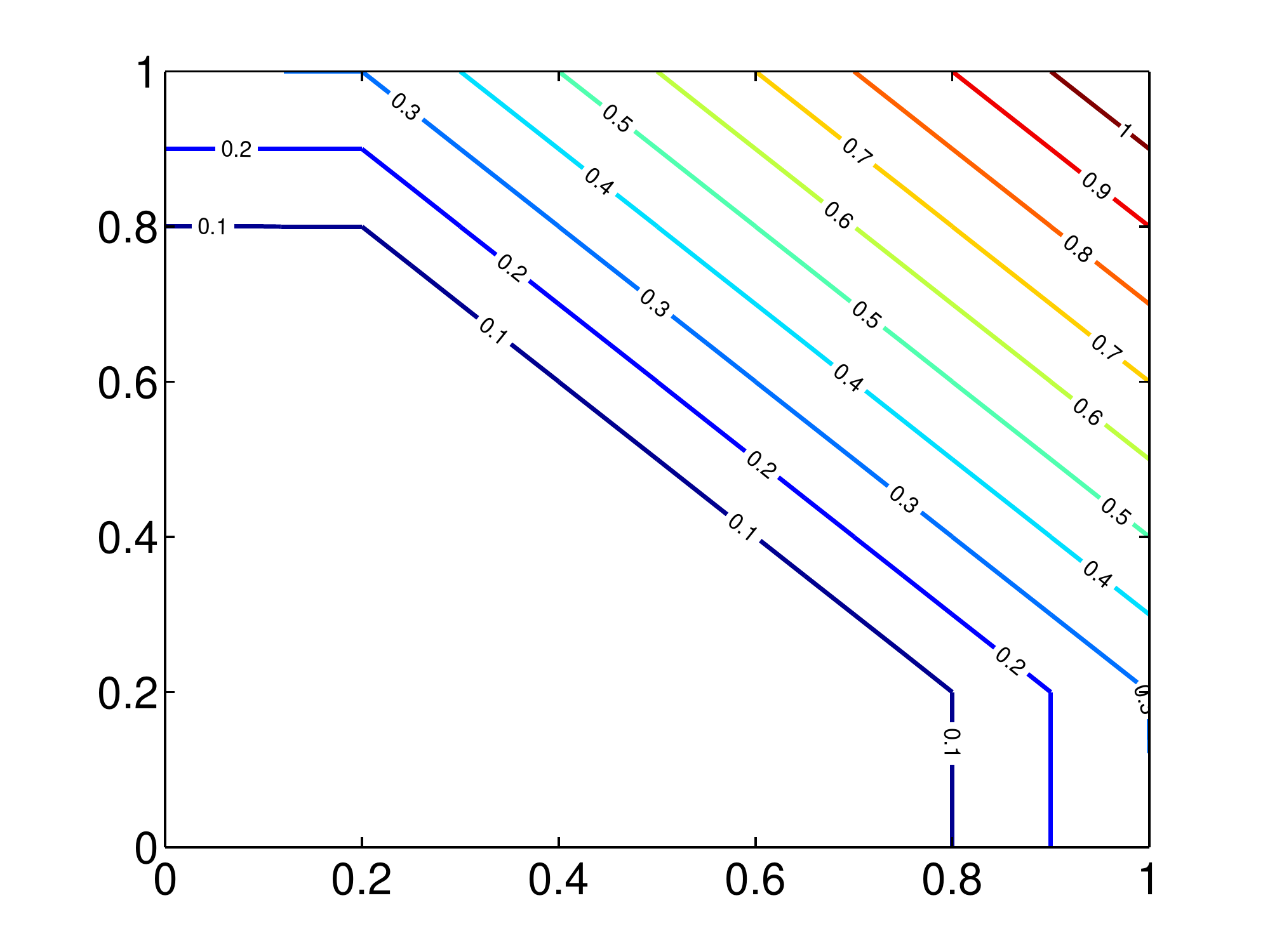}
\caption{Contour plot of the numerical solution of the Monopolist Problem, $n=21$.
}
\label{fig:mon1}
\end{figure}
	
\begin{table}[htdp]
\begin{center}
\begin{tabular}{|l | l | l | l | l | l | l |}
\hline 
 Method  \textbackslash ~n & 8 & 16 & 32 & 64 & 128 &  181
 \\  \hline
Zeroth Order Quadrature  & .14 & .071 & .010 & .0048 & .0042 &  .0055
\\  \hline
Trapezoidal Rule Quadrature & .05 & .005 & .01 & .005 & .0037 & .0008
\\  \hline
Method from~\cite{Morin1}  & .08 & .03 & .033 & .017 & x & x 
\\  \hline
\end{tabular} 
\end{center}
\caption{Monopolist Problem: error  ($L^\infty$) for the methods}

\begin{center}\begin{tabular}{| l | r | r | r | r | r | r |}
\hline 
Method  \textbackslash ~n & 8 & 16 & 32 & 64 & 128 &  181
 \\  \hline
Zeroth Order Quadrature &  & 2.1 & 1.7 & 7.4 & 69.3 & 417.7
\\  \hline
Trapezoidal Rule Quadrature & 2.4 & 2.2 & 1.3 & 5.6 & 80.2 & 455.1 
\\  \hline
Method from~\cite{Morin1}  &  0.2 & 1.0 & 17.0 & 751.0 & x & x
\\ \hline
\end{tabular} 
\caption{Monopolist Problem: run time (s) for the methods}

\end{center}
 \label{table:mon1}
\end{table}

\section{Reformulating polyhedral objective functions as linear constraints}
In this section we review a technique for  reformulating various convex optimization problems. This is needed to translate the optimization problem into standard forms (Quadratic and Linear Programs).  This reformulation is necessary for the (faster) MOSEK implementation.  On the other hand, the CVX implementation recognizes convex optimization problems, and does not require reformulation.

We also observed that different formulations of the optimization problems can lead to varying computational performance.

We start with a method for writing $L^\infty$ and $L^1$ norms (of either gradients or the variables) as a linear program (i.e. linear constraints and linear objectives).  
This is a standard technique see for example, see~\cite{BoydBook}
and~\cite{Morin1}.

\subsection{$\ell^\infty$ and $\ell^1$ projection}
The reformulation of the objective function 
$
\max_k\abs{u_k-g_k}
$
 is accomplished by setting the objective function to $t$ and enforcing,  $t = \max\abs{u-g}$, using the linear constraints
 \begin{align*}
g_k - t \le u_k &\le g_k + t,\quad \text{ for each $k$. }
\end{align*}

Similarly,  the reformulation of the objective function $\sum_k \abs{u_k-g_k}$ is accomplished by  setting the objective function to
$
 \sum_k t_k,
$
and enforcing $\abs{ u_k - g_k} = t_k$  
using the linear constraints 
\begin{align*}
g_k - t_k \le u_k &\le g_k + t_k, \quad \text{ for each $k$. }
\end{align*}

\subsection{Pointwise Constraints}
Dirichlet boundary conditions 
\[
u(x) = f(x), \quad x\in \partial D,
\]
can be implemented using equality constraints.

For variational problems with one degree of freedom, (e.g., $H^1$
projection) one equality constraint should be  included  to force uniqueness of the
minimizer, for example  by setting
\[
u(0) = 0.
\]
The latter constraint is not necessary for the problem to be well-posed, but it improves the conditioning of the problem and can speed up the solver.
For example, in the case of $H^1$ projection in one dimension, with $n =2001$, the  solution time improved from 32.3s to 8.7s.

\section{Performance improvement via Conic Programming}
Conventional wisdom is that quadratic programming is generally
faster than conic programming, and this is certainly the case for
the $L^2$ projections.  

However, using QP, the $H^1$ projections  were much slower than the
$L^2$ projections.  On the other hand, \cite{Morin1} used a conic
reformulation and found these solutions times of the $H^1$ and $L^2$ projections to be comparable (although in both cases slower than ours).  

The conic reformulation replaces the objective 
$
\frac 1 2 u (D_x^T D_x) u
$
with the new variable $t$ and the constraint
$
t \le  (D_x u)^2.
$
This formulation takes advantage of a matrix factorization 
\[
D_{xx} = D_x^T D_x.
\]
Our implementation of the QP did not take advantage of the matrix factorization, which may explain the performance difference.

We reformulated the quadratic objective as a conic
constraint.  Computed this way, the $H^1$ projection was no more costly than the $L^2$ projection. The results are shown in
Tables~\ref{tableL2}--\ref{tableH12d}.

\begin{table}[t]

\begin{tabular}{|c||c|c|c|c|}
\hline
n & 501 & 1001 & 2001 & 4001 
\\ \hline \hline
CP time & 1.4 & 1.4 & 1.5 & 1.8
\\ \hline
QP time & 1.4 & 1.5 & 1.7 & 1.8 
 \\ \hline
\end{tabular}
 \caption{Run time for $d=1$, $L^2$ projection, using Quadratic Program and Conic Program.}
 \label{tableL2}

\begin{tabular}{|c||c|c|c|c|}
\hline 
method / n & 501 & 1001 & 2001  & 4001 
\\ \hline
   \hline
CP time & 1.4 & 1.4 & 1.4  & 1.8
\\ \hline
QP time & 1.5 & 2.3 & 9.7  & 59
 \\ \hline
\end{tabular}
 \caption{Run time for $d=1$, $H^1$ projection, using Quadratic Program and Conic Program.}
 \label{tableH1}

\begin{tabular}{|c||c|c|c|c|c|c|c|c|c|}
\hline
n            & 4     &  8  & 16   & 20   & 32 & 45 &64 & 90 & 128
\\ \hline \hline
CP time  & 1.6 &1.7 & 1.9  & 2.3      & 5.5 & 13 & 42 &101 & 251
\\ \hline
QP time  & 1.4 & 1.5 & 2.3 & 6.4     & 117  & 1090  & x   & x &x
 \\ \hline
\end{tabular}
 \caption{Run time for $d=2$,  $H^1$ projection, using Quadratic Program and Conic Program.}
 \label{tableH12d}
\end{table}

\section{Conclusions}
In this article we considered the problem of approximating the solution
of variational problems subject to the constraint that the admissible
functions must be convex.   This problem has applications in shape optimization, and in mathematical economics.  It is related to the optimal mass transportation problem. 

This problem has proven to be computationally  challenging.  Counterexamples show that earlier approaches can fail to correctly represent convex functions, or require very costly techniques.

We introduced a polyhedral approximation of the cone of convex functions.  This approximation is
computationally efficient in the sense that it can be represented by a relatively small 
number of linear inequalities.  

We established an estimate which quantifies the approximation error (the degree of non-convexity or uniform convexity) based on the number of directions in which convexity was tested.   Numerical results showed that for solutions of problems which are not strictly convex, more directions are needed.  However in many cases where solutions are grid aligned or uniformly convex only a few directions are required. 

Our methods are computationally efficient: we have approximated the intractable convexity constraint by a relatively small number of linear constraints.     The increased efficiency is due to the reduction in the number of constraints to enforce (approximate) convexity.  The results are significantly (hundreds of times) faster, and generally more accurate compared to the most efficient method previously available.

Well resolved two dimensional  problems can be solved by easily obtained optimization packages.  In particular, we computed the  solution of the Monopolist problem with enough accuracy to resolve the gradient mapping.

\bibliographystyle{alpha}
\bibliography{ConvexVariational}
\end{document}